\documentclass[english,11pt]{article}

\usepackage[german,french,english]{babel}
\usepackage[cp850]{inputenc}
\usepackage{latexsym,graphicx, fancybox}
\usepackage{amssymb, amsmath, amsfonts}
\usepackage{color}
\usepackage{latexsym}

\font\tenmath=msbm10 scaled 1200
\font\sevenmath=msbm7 scaled 1200
\font\fivemath=msbm5 scaled 1200 

\newfam\mathfam \textfont\mathfam=\tenmath
\scriptfont\mathfam=\sevenmath \scriptscriptfont\mathfam=\fivemath

\def\R{{\mathbb R}}
\def\N{{\mathbb N}}
\def\E{{\mathbb E}}

\def\P{{\mathbb P}}

\def\Q{{\mathbb Q}}

\newtheorem{Theorem}{Theorem}[section]
\newtheorem{Definition}[Theorem]{Definition}
\newtheorem{Proposition}[Theorem]{Proposition}

\newtheorem{Lemma}[Theorem]{Lemma}
\newtheorem{Corollary}[Theorem]{Corollary}
\newtheorem{Remark}[Theorem]{Remark}
\newtheorem{Example}[Theorem]{Example}

\newfam\mathfam \textfont\mathfam=\tenmath
\scriptfont\mathfam=\sevenmath \scriptscriptfont\mathfam=\fivemath

\def \^#1{\if#1i{\accent"5E\i}\else{\accent"5E#1}\fi}

\def \ind {1 \mkern -5mu \hbox{I}}

\def \b{\beta}

\def \r{\rho}
\def \s{\sigma}

\begin{document}
\selectlanguage{english}
\title{\bf Convex ordering of solutions to one-dimensional SDEs}

\author{ 
{\sc Benjamin Jourdain} \thanks{CERMICS, Ecole des Ponts, INRIA, Marne-la-Vall\'ee, France. E-mail: {\tt   benjamin.jourdain@enpc.fr}}~\footnotemark[3]
\and   
{\sc  Gilles Pag\`es} \thanks{Laboratoire de Probabilit\'es, Statistique et Mod\'elisation, UMR~8001, Campus Pierre et Marie Curie, Sorbonne Universit\'e case 158, 4, pl. Jussieu, F-75252 Paris Cedex 5, France. E-mail: {\tt  gilles.pages@sorbonne-universite.fr}}~\thanks{This research
benefited from the support of the ``Chaire Risques Financiers'', Fondation du Risque}}
\date{}

\maketitle 
 \begin{abstract}
   In this paper, we are interested in the propagation of convexity by the strong solution to a one-dimensional Brownian stochastic differential equation with coefficients Lipschitz in the spatial variable uniformly in the time variable and in the convex ordering between the solutions of two such equations. We prove that while these properties hold without further assumptions for convex functions of the processes at one instant only, an assumption almost amounting to spatial convexity of the diffusion coefficient is needed for the extension to convex functions at two instants. Under this spatial convexity of the diffusion coefficients, the two properties even hold for convex functionals of the whole path. For directionally convex functionals, the spatial convexity of the diffusion coefficient is no longer needed.  Our method of proof consists in first  establishing the results for time discretization schemes of Euler type and then transferring them to their limiting  Brownian diffusions. We thus exhibit approximations which avoid {\em convexity arbitrages} by preserving convexity propagation and comparison and can be computed by Monte Carlo simulation.
 \end{abstract}
 \renewcommand{\abstractname}{Abstract}

\section{Introduction}
{In this paper, we carry on in-depth investigations about various aspects of propagation of convexity (and monotone convexity) by the flow of a one-dimensional Brownian diffusion
\begin{equation*}
dX_t=b(t,X_t)dt+\sigma(t,X_t)dW_t
\end{equation*} where $b,\sigma:[0,T]\times\R\to\R$ and $(W_t)_{t\ge 0}$ is a one-dimensional Brownian motion as well as the comparison  results with respect to marginal or functional convex orderings between two such diffusions. Studies on convex ordering(s) for stochastic processes go back to more than fifty years  with a seminal work by Hajek in~\cite{Hajek} (for increasing convex ordering), followed in the late 1990's by El Karoui, Jeanblanc and Shreve in~\cite{EKJS},  Martini in~\cite{Martini} and  then by various authors, among them Schied and Stadje in~\cite{SS}, Bergenthum and R\"uschendorf in a series of papers~\cite{Ruschendorf1, Ruschendorf2,Ruschendorf3} (on which we will comment further in Remark~\ref{rem:contrex}) or, more recently,~\cite{Pag2016} (for various stochastic dynamics),~\cite{LiuPagAAP,LiuPag} (for McKean-Vlasov processes) or~\cite{JourPagVolt22} (for Volterra processes). Except for~\cite{Hajek}, most of these papers were   motivated  by   financial applications to the pricing and hedging of derivative products in a complete market where the underlying traded asset price is supposed to be driven by/obey a local volatility model.
}

\bigskip
{In Section \ref{sec:funconv}, we state a theorem (Theorem~\ref{thm:start}) which combines the state of the art about convexity propagation and functional convex ordering(s) for strong solutions to $1$-dimensional Brownian SDEs combined with recent improvements developed while tackling similar problems for other classes of stochastic processes. This statement assumes convexity of at least one of the two diffusion coefficients in the spatial variable. The objective of this paper is to  generalize  this theorem essentially  in two directions. The first one is to  analyze the role played in dimension one by the convexity assumption made on the diffusion coefficient in the propagation of  convexity and subsequently in the comparison property. The second one is to extend the class of functionals for which such type of propagation/comparison can be established (see further on).   
}

{Our first objective is motivated by a result  from~\cite{EKJS}. In this paper, {\em without} assuming spatial convexity of the diffusion coefficient $\sigma$, the authors establish the propagation of convexity property (and a somewhat hidden comparison result through their proof of the {\em tracking error} property) for $1$-marginals of the diffusion~--~i.e. at the fixed time $T$~--  in the sense that for every convex function  $f:\R\to \R$
\[
x\longmapsto \E\, f(X^x_{_T})\quad \mbox{ is convex} ,
\]
where $(X^x_t)_{t\in[0,T]}$ denotes the solution to the above stochastic differential equation started from $X_0=x$.} In Section~\ref{subsec:CCsC},  we revisit this result in Theorem~\ref{thmmarg1} under weaker regularity assumptions on the coefficients 
  and establish the convexity of $x\mapsto \E\, f(X^x_{_T})$  under standard Lipschitz  assumptions only.  We prove accordingly under the same assumptions that the convex ordering property for the marginal distributions at time $T$ holds as well.

{To elucidate this paradox about the spatial convexity assumption needed on the diffusion coefficient, we were led to introduce the notion of $m$-marginal convex ordering that is comparison in convex order of the $m$-dimensional marginals of the diffusions.  By contrast with the purely marginal case $m=1$, we state in Section~\ref{subsec:Cm2} that, as soon as  
  convexity propagation holds for convex functions of the diffusion at $m$ instants with $m\ge 2$,   then $x\mapsto \E\, \sigma(t,X^x_t)$  is necessarily convex for every $t\!\in [0,T]$. In particular, in the autonomous case $\sigma(t,x)=\sigma(x)$, this implies that the function $\sigma$ itself is convex. No such spatial convexity is assumed by Schied and Stadje~\cite{SS} to extend the results in~\cite{EKJS} (convexity propagation and somewhat hidden comparison through robustness of the hedging strategy)  to payoffs depending on $m$ instants which are directionally convex. There is an overlap between the classes of directionally convex and convex functions but none contains the other, except in dimension one where they coincide. In our second main result (Theorem~\ref{thm:dirconv} in Section~\ref{sec:directconvex}), we check that propagation and comparison hold for directionally convex functionals.

{
On our way through the paper, we  always provide  results for both ``regular'' and  increasing convex propagation and ordering,
thus extending  the seminal theorem by Hajek~\cite{Hajek}.  
}

\smallskip 
{Our method of proof  for all these results always relies on the same paradigm: we first  establish  our results for time discretization schemes of Euler type and then transfer them using appropriate  convergence theorems of these schemes to their respective   Brownian diffusions. Like in~\cite{LiuPag}, but in a much more general and demanding case, we will need to consider Euler schemes with  {\em truncated noise} to obtain convexity propagation in discrete time. The truncation is chosen in such a way that the strong $L^1$-convergence to the limiting diffusion is preserved. Schemes involving time integrals of the coefficients converge without continuity of these coefficients in the time variable. Standard schemes avoiding such integrals are also considered under more stringent continuity assumptions. This is important in practice, especially in Finance when dealing with  derivative products (European options), since they provide approximations which preserve convexity propagation and comparison along convex functionals thus avoiding {\em convexity arbitrages} and can be computed by Monte Carlo simulation. 
  Also note that the results obtained in discrete time have their own interest for applications.
}

\medskip {The paper is organized as follows. Section \ref{sec:funconv} is devoted to convexity propagation for convex functionals and functional convexity comparison. In Section~\ref{sec:main}, we make a few definitions precise and  state the main results of the paper. In Theorem~\ref{thm:cvmarg} devoted to the $1$-marginal case, we show that spatial convexity of the diffusion coefficients can be relaxed. Then, in Proposition~\ref{prop:CNconv2}, we check that as soon as  the $2$-marginals functions come in the game, the assumptions to be made on the coefficients of the diffusions are nearly the same as those needed for general functionals of the whole paths of the processes.  The third part of this section is devoted to the class of directionally convex functionals which is not a subset of convex functionals (except in the $1$-marginal case) and for which  we establish  a theorem similar to Theorem~\ref{thm:start} without assuming spatial convexity of the diffusion coefficients. Section~\ref{sec:Eulertrunc} is devoted to several variants of truncated Euler schemes. They turn out to be the main tools used in our proofs since we first establish our results in discrete time and then transfer them to the continuous time setting by appropriate limit theorems. Section~\ref{sec:proofs} is devoted to the proofs of the main results stated in Section~\ref{sec:main}.}

\paragraph{Notation}
{
\begin{itemize}
\item We denote by ${\cal L}(X)$ the law of a random vector or a stochastic process $X$.
\item For $p\in[1,+\infty]$ (resp. $p\in(0,1)$), $\|X\|_p$ denotes the $L^p(\P)$-norm (resp. pseudo-norm) of the random vector $X$ defined on $(\Omega, {\cal A}, \P)$.
\item ${\cal C}([0,T], \R)$ denotes the space of continuous functions from $[0,T]$ to $\R$.
\item  We denote the sup-norm of a function $f\in{\cal C}([0,T], \R)$ by $\|f\|_{\infty}= \sup_{t\in [0,T]} |f(t)|\le +\infty$.
\end{itemize}
}
\section{Functional convex ordering of one-dimensional diffusions}\label{sec:funconv}
Before being more specific on the state  of the art functional convex ordering for one dimensional Brownian diffusions, let us briefly recall some definitions  about convex ordering of random vectors (in view of the proofs) and stochastic processes (which are main  objects of interest in our results).
\begin{Definition}[Convex orderings] Let  $U$, $V\!\in L^1(\P)$ be two $\R^m$-valued random vectors and let $X=(X_t)_{t\in [0,T]}$ and  $Y=(Y_t)_{t\in [0,T]}$ be  pathwise continuous real-valued stochastic processes such that $\E\,  \sup_{t\in [0,T]}|X_t|$ and $\E\,  \sup_{t\in [0,T]}|Y_t|<+\infty $.

\noindent {$(a)$ {\em Convex ordering}.    
We say that   $U$ is dominated for the convex ordering by $V$, denoted
\[
U\preceq_{cvx}V,
\]
 if, for every convex function $f:\R^m\to\R$, 
 \begin{equation}\label{eq:cvorder}
\E \, f(U) \leq \E\, f(V)\in \R\cup\{+\infty\}.
 \end{equation}
\smallskip
 \noindent $(b)$ {\em Monotone convex ordering}.  We say that   $U$ is dominated for the non-decreasing convex ordering by $V$, denoted
\[
U\preceq_{icv}V,
\]
if~\eqref{eq:cvorder} holds for every convex function $f:\R^m\to\R$ non-decreasing in each of its variables.
}

\smallskip
\noindent  {$(c)$ {\em Functional convex ordering}.    We say that  $X$ is dominated by  $Y$ for the functional  convex ordering, denoted 
$$
X\preceq_{cvx}Y
$$
  if for every  lower semi-continuous (l.s.c.) convex functional $F:{\cal C}([0,T], \R)\to \R$
\begin{equation}\label{eq:convF}
\E\, F(X)\le \E\, F(Y)\!\in \R\cup\{+\infty\}.
\end{equation}
$(d)$ {\em Functional monotone convex ordering}.We say that $X$ is dominated  by $Y$ (defined as above) for the increasing  functional  convex ordering, denoted 
$$
X\preceq_{icvx}Y,
$$ 
if~\eqref{eq:convF} holds true for every l.s.c. convex functional $F:{\cal C}([0,T], \R)\to \R$,  non-decreasing  for the natural partial order on ${\cal C}([0,T], \R)$~(i.e. for every $(x,y)\!\in {\cal C}([0,T], \R)\times {\cal C}([0,T], \R_+)$, $F(x)\le F(x+y)$).
}
\end{Definition}

Obviously, $U\preceq_{cvx}V\Longrightarrow U\preceq_{icvx}V$ and $X\preceq_{cvx}Y\Longrightarrow X\preceq_{icvx}Y$.

These definitions only involve the distributions of $U$, $V$ and $X$, $Y$ respectively. So we could have defined equivalently these orderings on the set of ${\cal P}_1(\R^m)$ and ${\cal P}_1\big(C([0,T], \R)\big)$ of distributions with finite first moment. 

Note that the expectations involved in the definitions make sense in $\R\cup\{+\infty\}$ since any convex function $f:\R^m\to\R$ (resp. l.s.c. convex functional $F:{\cal C}([0,T], \R)\to \R$) is bounded from below by an affine function (resp. by a continuous affine functional according to Lemma 7.5~\cite{AliBord}). The next proposition which can be proved by inf-convolution arguments (see the proof of~\cite[Lemma 1.9]{JourPagVolt22}) states a characterization of convex ordering(s)  based on the sole {\em Lipschitz} convex functions/functionals. 

\begin{Proposition}[{Convex ordering characterization}]\label{prop:ConLip} $(a)$ For every (non-decreasing in each of its variables) convex function $f:\R^m\to\R$, there exists a non-decreasing sequence of (non-decreasing in each of their variables) Lipschitz continuous convex functions having $f$ as a  pointwise limit.

\smallskip 
\noindent $(b)$  To establish (increasing) convex ordering  it suffices to establish~\eqref{eq:cvorder} for (non-decreasing in each of their variables) Lipschitz continuous convex  functions $f$.

\smallskip 
\noindent $(c)$
For every (non-decreasing) l.s.c. convex functional  $F:{\cal C}([0,T],\R)\to \R$ the exists a non-decreasing sequence of (non-decreasing) Lipschitz continuous convex functionals  having $F$ as a  pointwise limit.

\smallskip 
\noindent $(d)$  To establish (increasing) functional convex ordering  it suffices to establish~\eqref{eq:convF} for (non-decreasing) Lipschitz continuous convex  functionals $F$. \end{Proposition}

Let us now deal with the functional convex ordering of two one-dimensional diffusions
\begin{equation}\label{eq:SDE}
  X_t =X_ 0 +\int_0^t b(s,X_s)ds +\int_0^t \sigma(s,X_s)dW_s, \quad t\!\in [0,T],
\end{equation}
and 
\begin{equation}\label{eq:SDE2}
Y_t = Y_0 +\int_0^t \beta(s,Y_s)ds +\int_0^t \vartheta(s,Y_s)dW_s,\quad t\!\in [0,T], 
\end{equation}
where $b$, $\beta$, $\sigma$, $\vartheta : [0,T]\times \R\to \R$ are measurable functions, the initial positions  $(X_0,Y_0)$ and the  independent Brownian motion $W$ are defined on a probability space $(\Omega, {\cal A}, \P)$. We also assume that $b,\sigma,\beta,\vartheta$ satisfy the following hypothesis written for a generic function $f:[0,T]\times\R\to\R$:
\begin{align}\label{coeflip}
   {\rm Lip}(f):=\sup_{t\in[0,T],x\neq y\in\R}\frac{|f(t,x)-f(t,y)|}{|x-y|}<+\infty\mbox{ and }\sup_{t\in[0,T]}|f(t,0)|<+\infty,
\end{align}
which ensures that both diffusions have a unique adapted strong solution.

\medskip
Moreover (see e.g.~\cite[Proposition~7.2]{Pages2018} for a proof), if $X_0$ and $Y_0\!\in L^p(\P)$ for some $p>0$, then there exists a real constant $C= C_{b, \b,\s,\vartheta, p,T}$ only depending $b$, $\b$, $\s$,$\vartheta$, $p$ and $T$  such that
\begin{equation}\label{eq:Lpbound}
\Big\|\sup_{t\in [0,T]}|X_t \Big\|_p\le C (1+\|X_0\|_p) \quad \mbox{ and }\quad \Big\|\sup_{t\in [0,T]}|Y_t \Big\|_p\le C (1+\|Y_0\|_p).
\end{equation}

\medskip
\begin{Theorem}[Scaled martingale diffusions, see~\cite{PSP1},~\cite{LiuPagAAP},~\cite{JourPagVolt22}] \label{thm:start} {Assume that the functions $b$, $\beta$, $\sigma$, $\vartheta : \R\times [0,T]\to \R$ satisfy~\eqref{coeflip} and that $X_0,Y_0\in L^1(\P)$.}

\noindent {$(a)$ {\em Monotone convex ordering and convexity propagation}.  If \[
X_0\preceq_{icv} Y_0\quad \mbox{and}\quad \left\{\begin{array}{ll}
(i)_{\sigma,b} & \sigma(t,.):\R\to
\R_+\mbox{ and }b(t,.):\R\to\R\mbox{ are  convex for every $t\!\in [0,T]$},\\
\mbox{or}\\
(i)_{\vartheta,\beta} &\vartheta(t,.):\R\to\R_+\mbox{ and }\beta(t,.):\R\to\R \mbox{ are convex  for every $t\!\in [0,T]$},\\
\mbox{\and}\\
(ii) &\sigma(t,\cdot) \; \le\;
\vartheta(t,\cdot)\mbox{ and }b(t,.)\le\beta(t,.) \mbox{ for every $t\!\in [0,T]$,}
\end{array}\right.
\]
then
\begin {equation}\label{eq:comparXYm}
 X\preceq_{icv} Y.
\end{equation}
When $(i)_{\sigma,b}$ (resp. $(i)_{\vartheta,\beta}$) holds true then, for every l.s.c. convex $F : {\cal C}([0,T],\R) \to\R $  non-decreasing  for the natural partial order on   ${\cal C}([0,T],\R) $, 
\begin {equation}\label{eq:propagconvm}
\R\ni x\mapsto  \E\,F(X^{x})\in\R\cup\{+\infty\} (\mbox{resp. }\R\ni y\mapsto  \E\,F(X^{y})\in\R\cup\{+\infty\})\mbox{ is  convex and
non-decreasing},
\end{equation}
where $X^x$ (resp. $Y^y$) denotes the solution of~\eqref{eq:SDE} (resp.~\eqref{eq:SDE2}) starting from $X_0=x$ (resp. $Y_0=y$). 
}

\smallskip
\noindent   {$(b)$ {\em Convex ordering and convexity propagation for scaled  martingale diffusions}.  Assume that both SDEs share the  same drift affine in space of the form
\[
\forall\, t\!\in [0,T], \; \forall\, x\!\in \R,\quad b(t,x)=\beta(t,x)= \lambda(t)+\mu(t)x,
\]
where $\lambda,\, \mu:[0,T]\to \R$ are two bounded  Borel functions. 
}

{If   
\[
X_0\preceq_{cvx} Y_0\quad \mbox{and}\quad \left\{\begin{array}{ll}
(i)_{\sigma} & \sigma(t,.):\R\to
\R_+\mbox{ is  convex for every $t\!\in [0,T]$},\\
\mbox{or}\\
(i)_{\vartheta} &\vartheta(t,.):\R\to\R_+ \mbox{ is convex  for every $t\!\in [0,T]$},\\
\mbox{\and}\\
(ii) &\sigma(t,\cdot) \; \le\;
\vartheta(t,\cdot) \mbox{ for every $t\!\in [0,T]$,}
\end{array}\right.
\]
 then, 
\begin {equation}\label{eq:comparXY}
 X\preceq_{cvx} Y.
\end{equation}
}
\noindent When, $(i)_{\sigma}$ (resp. $(i)_{\vartheta}$) holds true, then, for every lower semi-continuous (l.s.c.) convex  functional  $F : {\cal C}([0,T],\R) \to\R $ on   ${\cal C}([0,T],\R)$, 
\begin {equation}\label{eq:propagconv}
\R\ni x\mapsto  \E\,F(X^{x})\in \R\cup\{+\infty\} (\mbox{resp. } \R\ni y\mapsto  \E\,F(X^{y})\in \R\cup\{+\infty\})\mbox{ is convex}.
\end{equation}
\end{Theorem}

We do not provide a detailed proof of this result in this paper since it appears as a synthesis of several more or less recently published papers. It combines the results established in~\cite{PSP1} (for true martingale diffusions as weak solutions of Brownian SDEs) and in~\cite{LiuPagAAP,LiuPag} (devoted to convex ordering for scaled martingale McKean-Vlasov diffusions and increasing  convex ordering for one-dimensional McKean-Vlasov diffusions respectively). Moreover, it benefits from various significant technical  improvements developed in~\cite{JourPagVolt22} (devoted to convex ordering for strong solutions of stochastic Volterra equations).

\bigskip
\noindent{\bf Remark.} The above theorem also has a {\em weak solution} version in which instead of assuming  that $b$, $\beta$, $\sigma$, $\vartheta$ satisfy~\eqref{coeflip} and $X_0,Y_0\in L^1(\P)$ one supposes that $b$, $\beta$, $\sigma$, $\vartheta$ are continuous on $[0,T]\times\R$ with at most affine growth in space uniformly in time and $X_0,Y_0\in L^p(\P)$ for some $p>1$.  
\section{Main results}\label{sec:main}
In order to state our main results, let us first introduce the $m$-marginal convex ordering between two stochastic processes.
\begin{Definition}[$m$-marginal convex ordering]
Let $m\!\in \N$.   A real-valued stochastic process $X=(X_t)_{t\in [0,T]}$ such that $\forall t\in[0,T],\;\E\,|X_t|<\infty$ is dominated by a real-valued process $Y=(Y_t)_{t\in [0,T]}$ such that $\forall t\in[0,T],\;\E\,|Y_t|<\infty$ for the $m$-marginal  convex ordering if, for every $t_1,\, t_2, \cdots,\,t_m\!\in [0,T]$ with $t_1\le t_2\le \cdots\le t_m$,
\[
(X_{t_1},\ldots,X_{t_m})\preceq_{cvx}(Y_{t_1},\ldots,Y_{t_m}).
\]
\end{Definition}

\subsection{$1$-marginal convex ordering and convexity propagation without  convexity of the diffusion coefficient}\label{subsec:CCsC}

For $1$-marginal convex ordering and convexity propagation, no spatial convexity property of the diffusion coefficients is needed.
\begin{Theorem}\label{thmmarg1}[1-marginal convex orderings without convexity of the diffusion coefficients] \label{thm:cvmarg} Assume that the functions $b$, $\beta$, $\sigma$, $\vartheta : \R\times [0,T]\to \R$ satisfy~\eqref{coeflip} and that $X_0,Y_0\in L^1(\P)$.
  
\smallskip
\noindent $(a)$ {\em Non-decreasing convex ordering}. Assume that either $\forall\, t\!\in [0,T]$, $b(t,\cdot)$ is convex or  $\forall\, t\!\in [0,T]$, $\beta(t,\cdot)$ is convex. We then have 
\[
\Big( b \le \beta \mbox{ and }0\le \sigma\le \vartheta \mbox{ on }[0,T]\times \R\mbox{ and }X_0\preceq_{icv} Y_0\Big)\Longrightarrow    X_{_T}\preceq_{icv}  Y_{_T}.
\]
Moreover, when $\forall\, t\!\in [0,T]$, $b(t,\cdot)$ is convex  (resp. $\forall\, t\!\in [0,T]$, $\beta(t,\cdot)$ is convex), then  for every non-decreasing convex function $f:\R\to \R$,
$$
x\mapsto \E\, f( X^{x}_{T})\in\R\cup\{+\infty\} \left(\mbox{resp. }y\mapsto \E\, f(Y^{y}_{T})\in\R\cup\{+\infty\}\right)  \; \mbox{  is non-decreasing and convex.}
$$

\noindent $(b)$ {\em Convex ordering}. Let  $b(t,x)=\beta(t,x) = \lambda(t)+\mu(t)x$, $(t,x)\!\in [0,T]\times \R$  be a time dependent affine function in space where $\lambda$ and $\mu$ are bounded Borel functions. 
Then, 

\[
\Big(0\le \sigma\le \vartheta\mbox{ on } [0,T]\times \R \mbox{ and } X_0\preceq_{cvx} Y_0\Big)\Longrightarrow  X_{_T}\preceq_{cvx} Y_{_T}.
\]
Moreover, for every convex function $f:\R\to \R$,
$$
x\mapsto \E\, f(X^{x}_{T})\in\R\cup\{+\infty\} \left(\mbox{resp. }y\mapsto \E\, f(Y^{y}_{T})\right)\in\R\cup\{+\infty\}  \; \mbox{  is convex.}
$$
\end{Theorem}

\bigskip
\noindent {\bf Remarks and comments.}  
$\bullet$ {\em From convex ordering to propagation of convexity and back.} The (non-decreasing) convexity of $x\mapsto \E\, f(X^x_{T})$ for $f:\R\to\R$ (non-decreasing) convex and the (increasing) convex ordering are closely intertwined. For the choice $(\beta,\vartheta)=(b,\sigma)$ and $X_0\sim\delta_{\alpha x+(1-\alpha )y}\le_{cvx} Y_0=\varepsilon \delta_x+(1-\varepsilon)\delta_y$ where $x,y\in \R$, $\alpha\in[0,1]$ and $\varepsilon$ independent of $W$ follows the Bernoulli distribution with parameter $\alpha$, the (increasing) convex  ordering implies that $\E \, f(X^{\alpha x +(1-\alpha)y}_{T})\le \alpha \E\, f(X^x_{T})+(1-\alpha) \E\, f(X^y_{T})$  for $f:\E\to\R$ (non-decreasing) convex while for the choice $X_0=x\le y=Y_0$ so that $X_0\le_{icv}Y_0$, the increasing convex  ordering implies that $\E\, f(X^x_{T})\le \E\, f(X^y_{T})$  for $f:\E\to\R$ non-decreasing and convex.
Conversely, we will rely on Euler discretizations of the SDEs to prove Theorem~\ref{thm:cvmarg} and at the level of these Euler schemes, the propagation of convexity turns out to be a key step to prove the convex ordering (see the proof of Proposition~\ref{prop:DirConvforEuler} below). 

\smallskip 
\noindent $\bullet$  {\em Theorem~\ref{thm:cvmarg} as an extension of Hajek's theorem.} As for the non-decreasing  convex ordering, the above comparison result extends a  theorem by Hajek~\cite{Hajek} in which one of the two diffusion coefficients  is supposed convex.
More precisely~\cite[Theorem 4.1]{Hajek} states that when 

\medskip
\begin{itemize}
\item $b,\sigma$ are measurable functions satisfying~\eqref{coeflip},

\smallskip
\item $\beta(t,x)$ and $\vartheta(t,x)$ do not depend on the time variable and are Lipschitz convex in the spatial variable,

\smallskip
\item $0\le \sigma(t,x)\le \vartheta(x)$
   and $b(t,x)\le \beta(x)$ for $(t,x)\in[0,T]\times\R$,
   
   \smallskip
   \item $X_0=x_0\le y_0=Y_0$,
\end{itemize}

\smallskip
\noindent then 
$$
X_T\preceq_{icv}Y_T.
$$

\medskip
\noindent $\bullet$ El Karoui, Jeanblanc and Shreve provide  in~\cite{EKJS}  a direct proof of the convexity of $x\mapsto \E\,f(X^x_{_T})$ for $f$ Lipschitz convex using stochastic calculus when $b(t,x)=r(t)x$, $\sigma(t,x)=x\gamma(t,x)$ with $\gamma:[0,T]\times (0,+\infty)\to\R$ continuous, bounded from above and such that $\partial_2\sigma(t,x)$ is continuous in $(t,x)$ and Lipschitz continuous and bounded in $x\in (0,+\infty)$ uniformly in $t\!\in[0,T]$. According to Schied and Stadje~\cite{SS} who relie on arguments based on partial differential equations, the convexity of $x\mapsto \E\,f(X^x_{_T})$ still holds for $f$ convex with polynomial growth when the drift vanishes and $\sigma(t,x)=x\gamma(t,x)$ is Lipschitz in $x$ uniformly in $t\in[0,T]$ with $\gamma:[0,T]\times (0,+\infty)\to (0,+\infty)$ continuously differentiable, bounded and bounded away from $0$. The regularity assumptions in the above theorem are weaker. For the reader's convenience, we extend in Appendix~\ref{App:A} the stochastic calculus argument in~\cite{EKJS} to prove the (non-decreasing) convexity of $x\mapsto \E\,f(X^x_{t_1},\cdots,X^x_{t_m})$ when $0\le t_1\le\cdots\le t_m\le T$ and $f:\R^m\to\R$ is directionally convex (and non-decreasing in each of its arguments). The case $m=1$ and $t_1=T$ correspond to the framework of the current section.}

\subsection{A necessary condition for $2$-marginal convex ordering}\label{subsec:Cm2} Let us now explain  why some stronger assumption on the diffusion coefficients is needed to get convex order at the level of two-dimensional marginals. An analogous result in discrete time is discussed in Remark \ref{rem:contrex}.

\begin{Proposition}[{$2$-marginal convex ordering}]\label{prop:CNconv2}
  $(a)$ Assume that $\sigma,\vartheta$ satisfy~\eqref{coeflip} and are continuous on $[0,T]\times\R$, that $b(t,x)=\beta(t,x)=\lambda(t)+\mu(t)x$ with continuous functions $\lambda,\mu:[0,T]\to\R$  and that $X_0,Y_0\in L^1(\P)$. 
If the solution $X=(X_t)_{t\in[0,T]}$ of~\eqref{eq:SDE} is dominated by the solution $Y=(Y_t)_{t\in[0,T]}$ of~\eqref{eq:SDE2} for the 2-marginal convex ordering, then $\E\,|\sigma(t,X_t)|\le \E\,|\vartheta(t,Y_t)|$ for each $t\in[0,T]$.
 
 \smallskip
 \noindent $(b)$ In particular, when $\vartheta=\sigma$, the above $2$-marginal convex ordering for all $X_0,Y_0\in L^1(\P)$ such that $X_0\le_{cvx} Y_0$ implies the convexity of $x\mapsto\E|\sigma(s,X^x_s)|$ for each $s\in[0,T]$.

 \smallskip
 \noindent {$(c)$  In the special case when  $\sigma$ is autonomous in the sense that $\sigma(t,x)=\sigma(x)$, the above $2$-marginal convex ordering implies that $|\sigma|:\R\to\R$ is a convex function.}
 \end{Proposition}

In fact, as shown by Schied et Stadje~\cite{SS} in the vector case, to establish a functional convex ordering and  propagation without such stronger assumption on the diffusion coefficient, we have to switch to a different class of functionals, namely directionally convex  functionals,  which induce  the directional convex ordering. This is the object of the next section.

\subsection{Ordering and propagation for directionally convex functionals}\label{sec:directconvex}
{Prior to the main result of this section (Theorem~\ref{thm:dirconv}), we first recall the definitions of  directional convexity  in both vector and functional settings and their first properties illustrated by some examples. In particular, we stress the fact that directional convexity is not a restriction 
  of the notion of convexity but a different class of functions/functionals (except in one dimension) with an overlap but none being more general than the other. We refer to Section 3.12 of~\cite{MulSto} where this notion is introduced and investigated for more details. As a second step, we define the resulting stochastic order, first for random vectors and then  for stochastic processes. Finally we state our main result. }
\begin{Definition}[Directional convexity]\label{def:directconv}

 \smallskip
 \noindent $(a)$ A function $f:\R^m\to\R$ is said to be directionally convex if it is Borel measurable and 
\begin{equation}\label{eq:directconvdef}
    \forall (x,y,z)\in\R^m\times\R_+^m\times\R_+^m,\quad f(x+y+z)-f(x+y)-f(x+z)+f(x)\ge 0.
\end{equation}

 \smallskip
 \noindent $(b)$ A functional $F:{\cal C}([0,T],\R)\to\R$ is said to be directionally convex if it is measurable and 
 \begin{equation}\label{eq:functconvdef}
 \forall (x,y,z)\in{\cal C}([0,T],\R)\times {\cal C}([0,T],\R_+)^2,\quad F(x+y+z)-F(x+y)-F(x+z)+F(x)\ge 0.
 \end{equation}

 \end{Definition}

\begin{Remark}[Characterisations]\label{remdircon}
   $\mathbf{P1.}$ A directionally convex function $f:\R^m\to\R$ is measurable and mid-point convex in each of its variables (the other variables being fixed) whence, according to ~\cite{Sierpinski1920}, convex in each of its variables. This implies that for each $i\in\{1,\cdots,m\}$, the right-hand derivative $\partial_{x_i+}f(x_1,\ldots,x_{i},\ldots,x_m)$ is non-decreasing in each of its variables and therefore locally bounded. As a consequence, the function $f$ is continuous. 

\smallskip
\noindent $\mathbf{P2.}$ A twice continuously differentiable function $f:\R^m\to\R$ is directionally convex if {and only if}  for all $i,j\in\{1,\cdots,m\}$, $\partial^2_{x_ix_j}f\ge 0$.

\smallskip
\noindent $\mathbf{P3.}$ When $d=1$, directional convexity amounts to regular convexity.

\smallskip
\noindent $\mathbf{P4.}$ When $g:\R\to\R$ is convex, then for each $i\in\{1,\cdots,m\}$, $\R^m\ni(x_1,\ldots,x_m)\mapsto g(x_i)$ is directionally convex.

\smallskip
\noindent {\noindent $\mathbf{P5.}$ A twice continuously Fr\'echet differentiable functional $F: {\cal C}([0,T], \R)\to \R$ is directionally convex if and only if, for every $x\!\in  {\cal C}([0,T], \R)$ and  every   $u,\,v\!\in {\cal C}([0,T], \R_+)$,
\[
D^2F(x).(u,v)\ge 0.
\]}
\end{Remark}

\noindent {The following examples illustrate that the notions of convexity and directional convexity have overlap but that none of the two is  more general than the other.}
\medskip
\noindent \begin{Example} \label{exple:dirconv}
  {$(a)$ The funtion $f:\R^2\to \R$ defined by $f(u,v)= |u-v|$ is (clearly) convex but not directionally convex. Indeed, if we set 
  \[
  x=y=(1,0),\quad z=(1,\tfrac 52)
  \]
 then $ f(x+y+z)-f(x+y)-f(x+z)+f(x)= |3- \frac 52 |-|2-0| -|2-\frac 52|+|1-0|= -1<0$ so that~\eqref{eq:directconvdef} is clearly violated.}

\medskip
\noindent {$(b)$ Let us consider the quadratic function $f:\R^2\to \R$ defined by
\[
f(u,v)= au^2+ 2\,cuv+bv^2\quad \mbox{where $a,b>0$}.
\]
One checks that  $f$ is convex in each of its variables and that
}
\begin{itemize}
\item {$f$ is convex if and only if $|c|\le \sqrt{ab}$,}

\smallskip
\item {$f$ is directionally convex if and only if $c\ge 0$ owing to the above Property {\bf P2}.}
\end{itemize}

\smallskip
{ Consequently, $f$ can be convex but not directionally convex (if $-\sqrt{ab}<c<0$), directionally convex without being convex (if $c>\sqrt{ab}$)  and both convex and directionally convex if $0\le c\le \sqrt{ab}$.}

\medskip
\noindent {$(c)$ By the same arguments,  this time based on  Property $\mathbf{P5}$, the functional $F: {\cal C}([0,T], \R)\to \R$ defined for every continuous function $x$ by 
\[
F(x)= \int_{[0,T]^2} (ax(s)^2 +2cx(s)x(t)+bx(t)^2)\,ds\,dt
\]
is
}
\begin{itemize}
\item  {convex if and only if $|c|\le \sqrt{ab}$,
}
\smallskip
\item {directionally convex if and only if $c\ge 0$.}
\end{itemize}
with the same conclusions as in~$(b)$.

\medskip 
\noindent {$(d)$ Let  $F: {\cal C}([0,T], \R)\to \R$ be the functional defined by 
\[
F(x) = \Psi\left( \int_0^T \varphi(x(s))ds\right), \quad x\!\in {\cal C}([0,T], \R),
\]
where $\Psi,\varphi:\R\to\R$  are two 
convex functions. }

\medskip
$\bullet$ {If $\Psi$  is non-decreasing then $F$ is convex,
}

\medskip
$\bullet$ {If both $\Psi$  and $\varphi$ are non-decreasing, then $F$ is directionally convex.} 

\medskip
{Consequently $F$ is ``more often'' convex than directionally convex. The proof, based on a direct checking of~\eqref{eq:functconvdef} is postponed to Appendix~\ref{app:B}.}      
 \end{Example}    
               The following lemma is the counterpart  for directionally convex functions of the  properties  of monotone approximation of convex functions by Lipschitz convex ones stated in Proposition \ref{prop:ConLip}. 
\begin{Lemma}\label{lem:approxdirconv}
   Let $f:\R^m\to\R$ be directionally convex (non-decreasing in each of its variables). Then there exists a non-decreasing sequence $(f_n)_{n\ge 1}$ of Lipschitz continuous directionally convex (non-decreasing in each variable) functions having $f$ as a pointwise limit. Moreover, for any $m$-dimensional random vector $U$ such that $U\in L^1(\P)$, $f^-(U)$ is integrable and $\E \,f(U)$ makes sense in $\R\cup\{+\infty\}$. 
 \end{Lemma}
 \noindent {\bf Proof.}
The directionally convex function $f$ being convex in each of its variables, it admits a left-hand derivative $\partial_{i-}f(u_1,\cdots,u_{i},\cdots,u_m)$ and a right-hand derivative $\partial_{i+}f(u_1,\cdots,u_{i},\cdots,u_m)$ with respect to its $i$-th coordinate for $i\in\{1,\cdots,m\}$. The directional convexity implies that $\partial_{i-}f$ and $\partial_{i+}f$ are non-decreasing in each of their variables and therefore locally bounded. As a consequence $f$ is locally Lipschitz continuous. The construction in~\cite[Theorem 3.12.7]{MulSto} then provides a non-decreasing sequence of Lipschitz continuous directionally convex functions $(f_n)_{n\ge 1}$ that converge pointwise to $f$ as $n\to+\infty$ : for $n\ge 1$, 
\begin{align*}
   f_n(u_1,\cdots&,u_m)=f(\psi_n(u_1,\cdots,u_m))\\&+\sum_{i=1}^m\left(\mathbf{1}_{\{u_i\le -n\}}\partial_{i-}f(\psi_n(u_1,\cdots,u_m))(u_i+n)+\mathbf{1}_{\{u_i\ge n\}}\partial_{i+}f(\psi_n(u_1,\cdots,u_m))(u_i-n)\right),
\end{align*} 
where $\psi_n(u_1,\cdots,u_m)=(-n\vee u_1\wedge n,\cdots,-n\vee u_m\wedge n)$. When $f$ is moreover non-decreasing in each of its variables, then so is $f_n$ since $\psi_n$ also is and $\partial_{i-}f$ and $\partial_{i+}f$ are non-negative.
                                        
\smallskip      
\noindent When $U\in L^1(\P)$, then, since $f_1$ is Lipschitz continuous, $f_1(U)$ is integrable. Since $f_1\le f$, $f^-\le f_1^-$ and we deduce that $f^-(U)$ is integrable.
\hfill$\Box$

An extension of this result to l.s.c. directionally convex functionals remains an open question to us. Moreover, we do not whether a l.s.c. directionally convex functional is bounded from below by a continuous affine functional. In particular, we do not know whether $\E\,F(X)$ makes sense in $\R\cup\{+\infty\}$ when $F:{\cal C}([0,T],\R)\to\R$ is l.s.c. directionally convex and $X=(X_t)_{t\in[0,T]}$ is a continuous real-valued process such that $\E\, \sup_{t\in [0,T]} |X_t|<+\infty $. That is why we restrict the class of directionally convex functionals in the next definition of directionally convex orderings.

  \begin{Definition}[Directionally convex orderings]

\noindent $(a)$ 
Let  $U$, $V\!\in L^1(\P)$ be two $\R^m$-valued random vectors.

\smallskip
\noindent {$\rhd$ Directionally convex ordering on $\R^m$}.   We say that   $U$ is dominated for the directionally convex ordering by $V$, denoted
\begin{equation}\label{eq:dirconvdef}
U\preceq_{dircvx}V
\end{equation}
 if~\eqref{eq:cvorder} holds for every directionally convex function $f:\R^m\to\R$.

\smallskip
\noindent {$\rhd$ Monotone directionally convex ordering on $\R^m$}.  We say that   $U$ is dominated for the non-decreasing directionally convex ordering by $V$, denoted
\[
U\preceq_{diricv}V
\]
 if~\eqref{eq:cvorder} holds for every directionally convex function $f:\R^m\to\R$ non-decreasing in each of its variables.

\smallskip
\noindent {$(b)$ Let $X=(X_t)_{t\in [0,T]}$  and $Y= (Y_t)_{t\in [0,T]}$ be two pathwise continuous  $\R$-valued stochastic processes satisfying $\E\, \sup_{t\in [0,T]} |X_t|^p +\E\, \sup_{t\in [0,T]} |Y_t|^p <+\infty $.} 

\smallskip
\noindent {$\rhd$ Functional directionally convex ordering}. We say that $X$ is dominated by $Y$  for the functional directionally convex ordering with growth index $p$, denoted 
$$
X\preceq_{dircvx(p)}Y,
$$  if for every continuous directionally convex functional $F:{\cal C}([0,T], \R)\to \R$  with polynomial growth of order $p$ (w.r.t. the sup norm)~(\footnote{i.e. $\exists C_{_F}<\infty$ such that for every $x\!\in {\cal C}([0,T], \R)\to \R$, $|F(x)|\le C_{_F}(1+\|x\|^p_{\sup})$.})  
\begin{equation}\label{eq:funcdirconvdef}
\E\, F(X)\le \E\, F(Y).
\end{equation}
\noindent {$\rhd$ Monotone functional directionally convex ordering}  We say that $X$ is dominated by $Y$  for the functional  increasing directionally   convex ordering with growth index $p$, denoted 
$$
X\preceq_{diricv(p)}Y,
$$  
if~\eqref{eq:funcdirconvdef} holds  for every non-decreasing continuous directionally convex functional $F:{\cal C}([0,T], \R)\to \R$  with polynomial growth of order $p$.
 \end{Definition}

\begin{Remark}\label{remdircvcv}
Property {\bf P4} in Remark~\ref{remdircon} implies that
 \begin{align}
   &(U_1,\cdots,U_m)\le_{dircvx}(V_1,\cdots,V_m)\Longrightarrow \forall i\in\{1,\cdots,m\},\; U_i\le_{cvx} V_i\label{dirmarg},\\
   &(U_1,\cdots,U_m)\le_{diricv}(V_1,\cdots,V_m)\Longrightarrow \forall i\in\{1,\cdots,m\},\; U_i\le_{icv} V_i\label{dirimarg}.
 \end{align}
\end{Remark}
We are now ready to state our result concerning directionally convex functional ordering.

\begin{Theorem}\label{thm:dirconv}
  Assume that the functions $b$, $\beta$, $\sigma$, $\vartheta : \R\times [0,T]\to \R$ satisfy~\eqref{coeflip} and that $X_0,Y_0\in L^p(\P)$ for some $p\in[1,+\infty)$.
  
\noindent $(a)$ {\em Non-decreasing directionally convex functional ordering}. 
Assume that either $\forall\, t\!\in [0,T]$, $b(t,\cdot)$ is convex or  $\forall\, t\!\in [0,T]$, $\beta(t,\cdot)$ is convex. Then,
\[
\Big(b\le\beta\mbox{ and }0\le \sigma\le \vartheta\mbox{ on } [0,T]\times \R\mbox{ and }X_0\preceq_{icv} Y_0\Big)\Longrightarrow X\preceq_{dircvx(p)} Y.
\]
Moreover, when $\forall\, t\!\in [0,T]$, $b(t,\cdot)$ is convex, then for each continuous non-decreasing directionally convex $F:{\cal C}([0,T],\R)\to\R$ with polynomial growth, $\R\ni x\mapsto \E F(X^x)\in\R$ is convex and non-decreasing.

\smallskip
\noindent $(b)$ {\em Directionally convex functional ordering}. Let  $b(t,x)=\beta(t,x) = \lambda(t)+\mu(t)x$, $(t,x)\!\in [0,T]\times \R$  be a time dependent affine function in space where $\lambda$ and $\mu$ are bounded Borel functions. 
Then,
\[
\Big(0\le \sigma\le \vartheta\mbox{ on } [0,T]\times \R\mbox{ and }X_0\preceq_{cvx} Y_0\Big)\Longrightarrow X\preceq_{diricv(p)}Y.
\]{Moreover, for each continuous directionally convex $F:{\cal C}([0,T],\R)\to\R$ with polynomial growth, $\R\ni x\mapsto \E F(X^x)\in\R$ is convex.}\end{Theorem}

\noindent {\bf Remark.}  According to Schied and Stadje~\cite{SS} who relie on arguments based on partial differential equations, $x\mapsto \E\,f(X^x_{t_1},\cdots,X^x_{t_m})$ is convex for $0\le t_1\le\cdots\le t_m\le T$ and $f:\R^m\to\R$ directionally convex with polynomial growth when the drift vanishes and $\sigma(t,x)=x\gamma(t,x)$ is Lipschitz in $x$ uniformly in $t\in[0,T]$ with $\gamma:[0,T]\times (0,+\infty)\to (0,+\infty)$ continuously differentiable, bounded and bounded away from $0$. The convexity of $x\mapsto \E\,f(X^x_{t_1},\cdots,X^x_{t_m})$ is derived in Appendix~\ref{App:A} by a generalization of the stochastic calculus argument in~\cite{EKJS} which requires more regularity on the coefficients of the SDE.

\section{Convergence and ordering of Euler schemes with truncated noise}\label{sec:Eulertrunc}

{Our proofs rely on the same ``transfer'' paradigm: prove the results  in discrete time on  discretization schemes   of  Euler type and then transfer them by a strong (or functional weak ) convergence theorem toward the diffusion. We will extensively use  for the questions under consideration in this paper some specific truncated Euler schemes, especially when dealing with monotone convex ordering and propagation. We  always consider a taylored Euler scheme  for which the lightest  assumptions are needed but, when it turns out not to be simulable, we also give a result for  a simulable Euler scheme under slightly more stringent assumptions.}

It turns out that we will need to truncate the Brownian increments in the Euler scheme to ensure propagation of convexity by the transitions associated with this scheme. That is why, in the next section, we address the convergence of such Euler schemes with truncated noise.
\subsection{Convergence of Euler schemes with truncated noise}
The (continuous  time) Euler scheme with step $\frac T m,\;m\ge 1$ that we use to discretize~\eqref{eq:SDE} is defined by freezing the spatial coefficients (and not the time coefficients) using $\underline t = \frac{kT}{m}$ when $t\!\in [ \frac{kT}{m},  \frac{(k+1)T}{m})$, $k=0,\ldots,m-1$ and $\underline T = T$ :
\begin{equation}\label{eq:Euler}
\bar X_t^m = X_0 +\int_0^t b(s, \bar X^m_{\underline s})ds +\int_0^t \sigma(s, \bar X^m_{\underline s})dW_s.
\end{equation}
If we denote by $t_k= \frac{kT}{m}$ the discretization times, we have $\bar X^m_0= X_0$ and 
\begin{equation}\label{eq:Eulerdisc}
\bar X^m_{t_{k+1}}= \bar X^m_{t_k} + \int_{t_k}^{t_{k+1}}b(s, \bar X^m_{t_k})ds + \int_{t_k}^{t_{k+1}}\sigma(s, \bar X^m_{t_k}) dW_s, \; k=0,\ldots,m-1.
\end{equation}
Introducing the independent standard Gaussian random variables $G_{k+1} = (W_{t_{k+1}}-W_{t_{k}})/ \sqrt{T/m}$, $k=0,\ldots,m-1$, we can also associate a truncated  Euler scheme relying on the truncated Brownian increments
\begin{equation}\label{eq:truncscheme}
Z^m_{k+1}= G_{k+1} \mathbf{1}_{\{|G_{k+1}| \le {\mathbf s}_m\}}, \quad k=0,\ldots,m-1
\end{equation}
i.e. $\widetilde X^m_0= X_0$ and 
\begin{equation}\label{eq:EulerTrunc}
\widetilde X^m_{t_{k+1}}= \widetilde  X^m_{t_k} + \int_{t_k}^{t_{k+1}}b(s, \widetilde  X^m_{t_k})ds + \sqrt{\int_{t_k}^{t_{k+1}}\sigma^2(s, \widetilde   X^m_{t_k})ds}\,Z^m_{k+1}, \; k=0,\ldots,m-1.
\end{equation}

\begin{Proposition}\label{prop:ConvfortEuler}
 Assume $b,\sigma$ satisfy~\eqref{coeflip} and let $p\in[1,+\infty)$. There exists a finite constant $C_p$ such that for each initial random variable $X_0$  and each $m\ge 1$, the truncated Euler scheme~\eqref{eq:EulerTrunc} satisfies
$$
{\cal W}_p({\cal L}(X_{t_0},\ldots,X_{t_m}),{\cal L}(\widetilde X^m_{t_0},\ldots,\widetilde X^m_{t_m}))\le C_p\left(m^{-\frac12}+{\Big(\frac{m}{{\mathbf s}_m}e^{-\frac{{\mathbf s}_m^2}{2}}\Big)^{\frac{1}{2p}}}\right)(1+\|X_0\|_p),
$$
  where, in the definition of the Wasserstein distance ${\cal W}_p$, $\R^{m+1}$ is endowed with $\displaystyle \|(x_0,\ldots,x_{m})\|_\infty=\max_{k=0,\ldots,m}|x_k|$.
\end{Proposition}
\begin{Remark}
   We deduce that, as $m\to+\infty$,  
   $$
   \liminf_m\frac{{\mathbf s}_m}{\sqrt{\ln m}}>\sqrt{2}\Rightarrow {\cal W}_p({\cal L}(X_{t_0},\ldots,X_{t_m}),{\cal L}(\widetilde X^m_{t_0},\ldots,\widetilde X^m_{t_m}))\to 0.
   $$
\end{Remark}

\medskip
\noindent {\bf Proof of Proposition~\ref{prop:ConvfortEuler}.}
To analyse the convergence of the Euler scheme with truncated Brownian increments, it is convenient to introduce a third scheme $(\hat X^m_{t_k})_{k=0,\ldots,m}$ inductively defined by $\hat X^m_{0}=X_0$ and 
$$\hat X^m_{t_{k+1}}= \hat X^m_{t_k} + \int_{t_k}^{t_{k+1}}b(s, \hat X^m_{t_k})ds + \sqrt{\int_{t_k}^{t_{k+1}}\sigma^2(s, \hat X^m_{t_k})ds}\,G_{k+1}, \; k=0,\ldots,m-1.$$
Since $\bar X^m_{t_0}=\hat X^m_{t_0}=X_0$ and the conditional law of $\bar X^m_{t_{k+1}}$ given $(\bar X^m_{t_0},\ldots,\bar X^m_{t_k})$ is Gaussian with expectation $\int_{t_k}^{t_{k+1}}b(s, \bar X^m_{t_k})ds$ and variance $\int_{t_k}^{t_{k+1}}\sigma^2(s, \bar X^m_{t_k})ds$ and the same statement holds with $\bar X^m$ replaced by $\hat X^m$, we have {by induction}  that ${\cal L}(\bar X^m_{t_0},\bar X^m_{t_1},\ldots,\bar X^m_{t_m})={\cal L}(\hat X^m_{t_0},\hat X^m_{t_1},\ldots,\hat X^m_{t_m})$.

By standard arguments, the fact that $b,\sigma$ satisfy~\eqref{coeflip} implies the existence of a finite constant $\kappa^{(1)}_{b,\sigma,p,T}$ such that
\begin{align*}
  &\Big \|\sup_{t\in [0,T]} |X_t|  \Big\|_p+ \Big \|\sup_{t\in [0,T]}|\bar X^m_t|  \Big\|_p+  \Big \|\max_{k=0,\ldots,m}|\hat X^m_{t_k}|  \Big\|_p +  \Big \|\max_{k=0,\ldots,m}|\widetilde X^m_{t_k}|  \Big\|_p\le \kappa^{(1)}_{b,\sigma,p,T}(1+\|X_0\|_p)\notag\\
  &\mbox{ and }\Big \|\sup_{t\in [0,T]} |X_t-\bar X^m_t|  \Big\|_p\le \frac{\kappa^{(1)}_{b,\sigma,p,T}}{\sqrt{m}}(1+\|X_0\|_p).
\end{align*}
Since \begin{align*}
  {\cal W}_p({\cal L}(X_{t_0},\ldots,X_{t_m}),{\cal L}(\widetilde X^m_{t_0},\ldots,\widetilde X^m_{t_m}))&\le {\cal W}_p({\cal L}(X_{t_0},\ldots,X_{t_m}),{\cal L}(\bar X^m_{t_0},\ldots,\bar X^m_{t_m}))\\&\phantom{\le}+{\cal W}_p({\cal L}(\hat X_{t_0},\ldots,\hat X_{t_m}),{\cal L}(\widetilde X^m_{t_0},\ldots,\widetilde X^m_{t_m}))\notag\\
  &\le \Big \|\sup_{t\in [0,T]} |X_t-\bar X^m_t|  \Big\|_p+\Big \|\max_{k=0,\ldots,m}|\hat X^m_{t_k} - \widetilde X^m_{t_k}|  \Big\|_p,
\end{align*}
to conclude, it is enough to check that
\begin{equation}
   \Big \|\max_{k=0,\ldots,m}|\hat X^m_{t_k} - \widetilde X^m_{t_k}|  \Big\|_p\le 2\kappa^{(1)}_{b,\sigma,2p}{\Big(\frac{m}{{\mathbf s}_m}e^{-\frac{{\mathbf s}_m^2}{2}}\Big)^{\frac{1}{2p}}}(1+\|X_0\|_p),  
\label{eultrunceul}\end{equation} which we now do.
To compare $(\hat X^m_{t_k})_{k=0,\ldots,m} $ and $(\widetilde X^m_{t_k})_{k=0,\ldots,m}$ we simply note that 
\[
(\hat  X^m_{t_k})_{k=0,\ldots,m}= (\widetilde X^m_{t_k})_{k=0,\ldots,m} \quad \mbox{ on the event  }\quad \Omega_m = \bigcap_{k=1}^m \left\{ |G_k|\le  {\mathbf s}_m\right\}.
\]
Then, when $\|X_0\|_{2p}<\infty$, 
\begin{align*}
  \Big \|\max_{k=0,\ldots,m}|\hat X^m_{t_k} - \widetilde X^m_{t_k}|  \Big\|_p &=  \Big \|\mathbf{1}_{\Omega^c_m} \max_{k=0,\ldots,m}|\hat X^m_{t_k} - \widetilde X^m_{t_k}|  \Big\|_p\\
  &\le \P(\Omega^c_m)^{\frac{1}{2p}} \Big \| \max_{k=0,\ldots,m}|\hat X^m_{t_k}|+ \max_{k=0,\ldots,m}|\widetilde X^m_{t_k}|  \Big\|_{2p}\\
  & \le \P(\Omega^c_m)^{\frac{1}{2p}}\left(\Big \| \max_{k=0,\ldots,m}|\hat X^m_{t_k}|   \Big\|_{2p}+\Big \| \max_{k=0,\ldots,m}|\widetilde X^m_{t_k}|  \Big\|_{2p}\right)\\
  &  \le 2\kappa^{(1)}_{b,\sigma,2p}  \P(\Omega^c_m)^{\frac{1}{2p}}\big(1+\|X_0\|_{2p}\big). 
\end{align*}
In particular, denoting by $\hat X^{m,x_0}_{t_k}$ and $\widetilde X^{m,x_0}_{t_k}$ the regular and the truncated Euler schemes both {starting}  from the deterministic initial position $x_0\in\R$, we have
$$\E\left[\max_{k=0,\ldots,m}\left|\hat X^{m,x_0}_{t_k} - \widetilde X^{m,x_0}_{t_k}\right|^p\right]  \le (2\kappa^{(1)}_{b,\sigma,2p})^{p}  \P(\Omega^c_m)^{\frac{1}{2}}\big(1+|x_0|\big)^p.$$
Integrating this inequality in $x_0$ with respect to the distribution of $X_0$ and using the independence of this random variable and the vector of Brownian increments $(G_1,\ldots,G_m)$, we deduce that
$$ \Big \|\max_{k=0,\ldots,m}|\hat X^m_{t_k} - \widetilde X^m_{t_k}|  \Big\|_p\le 2\kappa^{(1)}_{b,\sigma,2p} \P(\Omega^c_m)^{\frac{1}{2p}}\|1+|X_0|\|_p\le 2\kappa^{(1)}_{b,\sigma,2p} \P(\Omega^c_m)^{\frac{1}{2p}}(1+\|X_0\|_p).$$
Since, by the elementary inequality $\P(|G_1|> z) = 2\P(G_1> z) \le \left(1\wedge \sqrt{\frac{2}{\pi z^2}}\right)e^{-\frac{z^2}{2}}$, $z\ge 0$,
\begin{align*}
  \P(\Omega^c_m)  = \P \left(\max_{k=1,\ldots,m} |G_k |> {\mathbf s}_m\right) \le m\P \left(|G_1|> {\mathbf s}_m\right)
  	& \le \sqrt{\frac{2}{\pi}}\frac{m}{{\mathbf s}_m}  e^{-\frac{{\mathbf s}_m^2}{2}}, 
\end{align*}
we deduce~\eqref{eultrunceul}. \hfill$\Box$

\bigskip
To transfer our discrete time results on convex ordering to continuous time, we use the ${\cal C}([0,T], \R)$-valued interpolation operator $i_m$ associated to the mesh $(t_k = \frac{kT}{m})_{k=0,\ldots,m}$ already introduced in~\cite{Pag2016}
\[
i_m: \R^{m+1}\ni (x_{0},\ldots,x_m) \longmapsto \Big(t \mapsto  \sum_{k=1}^m \mathbf{1}_{[t_{k-1},t_k)}(t) \left(\frac{t_{k}-t}{t_{k}-t_{k-1}}x_{k-1} + \frac{t-t_{k-1}}{t_{k}-t_{k-1}} x_k\right)+\mathbf{1}_{\{t=t_n\}}x_m\Big)\!.
\]

\begin{Corollary}\label{coro:cvinterpol}Let us assume that $b,\sigma$ satisfy~\eqref{coeflip} and that $\|X_0\|_p< +\infty$ for some $p\in[1,+\infty)$. Then the truncated Euler scheme~\eqref{eq:EulerTrunc} $(\widetilde X^m_{t_k})_{k=0,\ldots,m}$ satisfies 
  \[
 {\cal W}_p\left({\cal L}(X),{\cal L}\left(i_m\big((\widetilde X^m_{t_k})_{k=0,\ldots,m}\big)\right)\right)\to 0\quad \mbox{ as }\quad m\to +\infty\mbox{ when }{\liminf_m\frac{{\mathbf s}_m}{\sqrt{\ln m}}>\sqrt{2}}.
\]   
\end{Corollary}
\noindent {\bf Proof.} For $(x_k)_{k=0,\ldots,m}\in(\R^{d})^{m+1}$, the function $i_m((x_k)_{k=0,\ldots,m})$ being continuous and piecewise affine with affinity breaks at times $t_k$ where it is equal to $x_k$, \begin{equation*}\| i_m((x_k)_{k=0,\ldots,m})\|_{\sup}= \max_{k=0,\ldots,m} |x_k|.
\end{equation*}
For $(x_k)_{k=0,\ldots,m},(y_k)_{k=0,\ldots,m}\in\R^{m+1}$, since $i_m((x_k)_{k=0,\ldots,m})-i_m((y_k)_{k=0,\ldots,m})=i_m((x_k)_{k=0,\ldots,m}-(y_k)_{k=0,\ldots,m})$, we deduce that $i_m$ is Lipschitz with constant $1$ from $\R^{m+1}$ to $\big({\cal C}([0,T], \R), \|\cdot\|_{\sup}\big)$ :
\[
\| i_m((x_k)_{k=0,\ldots,m}) -i_m((y_k)_{k=0,\ldots,m}) \|_{\sup}= \max_{k=0,\ldots,m} |x_k-y_k|.
\]
By integrating this equality in $((x_k)_{k=0,\ldots,m},(y_k)_{k=0,\ldots,m})$ against a ${\cal W}_p$ optimal coupling between ${\cal L}(X_{t_0},\ldots,X_{t_m})$ and ${\cal L}(\widetilde X^m_{t_0},\ldots,\widetilde X^m_{t_m})$ and using Proposition~\ref{prop:ConvfortEuler}, we deduce that
\[
 {\cal W}_p\left({\cal L}\left(i_m\big((X_{t_k})_{k=0,\ldots,m}\big)\right),{\cal L}\left(i_m\big((\widetilde X^m_{t_k})_{k=0,\ldots,m}\big)\right)\right)\le C_p\left(m^{-\frac12}+{\Big(\frac{m}{{\mathbf s}_m}e^{-\frac{{\mathbf s}_m^2}{2}}\Big)^{\frac{1}{2p}}}\right)(1+\|X_0\|_p),
 \] 
 where the right-hand side goes to $0$ as $m\to+\infty$ when $\liminf_m\frac{{\mathbf s}_m}{\sqrt{\ln m}}>\sqrt{2}$.
With the triangle inequality for ${\cal W}_p$ and the inequality ${\cal W}_p\left({\cal L}(X),{\cal L}\left(i_m\big((X_{t_k})_{k=0,\ldots,m}\big)\right)\right)\le \Big\|  \big \| X-i_m\big( (X_{t_k})_{k=0,\ldots, m}\big)\big\|_{\sup}\Big\|_p$ we conclude that it is enough to check that \begin{equation}
   \Big\|  \big \| X-i_m\big( (X_{t_k})_{k=0,\ldots, m}\big)\big\|_{\sup}\Big\|_p\to 0\quad\mbox{as}\quad m\to+\infty\label{convinterp}
\end{equation} to complete the proof.
The uniform continuity modulus $w(\xi, \delta)=\sup_{0\le s\le t\le (s+\delta)\wedge T}|\xi(t)-\xi(s)|$ of a function $\xi\!\in {\cal C}([0,T], \R)$ converges to $0$ as $\delta\to 0$. One easily checks that 
 \[
 \big \| X-i_m\big( (X_{t^m_k})_{k=0,\ldots, m}\big)\big\|_{\sup} \le  w\big(X, \tfrac Tm\big)\le 2\big \| X\big\|_{\sup}.
 \]
As, by standard estimations, $\big\| \|X\|_{\sup}\big\|_p \le C(1+\|X_0\|_p)<+\infty$,~\eqref{convinterp} holds by dominated convergence. \hfill$\Box$

\bigskip
No continuity of the coefficients $b$, $\sigma$ in their time variable is required for the law of the Euler scheme~\eqref{eq:EulerTrunc} with truncated Brownian increments to converge to {that} of the solution of the SDE~\eqref{eq:SDE}. The counterpart is that one needs closed-form formulas for the integrals $\int_{t_k}^{t_{k+1}}b(s,x)ds$ and $\int_{t_k}^{t_{k+1}}\sigma^2(s,x)ds$ to implement this scheme in practice. If such formulas are not available, one should rather rely on the {the  standard  Euler scheme} where the time variable of the coefficients is also discretized : $\widetilde X^m_0= X_0$ and 
\begin{equation}\label{eq:EulerTrunc2}
\widetilde X^m_{t_{k+1}}= \widetilde  X^m_{t_k} +\frac{T}{m} b(t_k, \widetilde  X^m_{t_k}) + \sqrt{\frac{T}{m}}\sigma(t_k, \widetilde   X^m_{t_k})\,Z^m_{k+1}, \; k=0,\ldots,m-1.
\end{equation}
The corresponding Euler scheme with non-truncated increments {reads}
 \begin{equation}\label{eq:Euler2}
\bar X_t^m = X_0 +\int_0^t b({\underline s}, \bar X^m_{\underline s})ds +\int_0^t \sigma({\underline s}, \bar X^m_{\underline s})dW_s
\end{equation}
and evolves inductively from one discretization time to the next according to 
\begin{equation}\label{eq:Eulerdiscdisc}
\bar X^m_{t_{k+1}}= \bar X^m_{t_k} +\frac{T}{m} b(t_k, \bar X^m_{t_k})+ \sqrt{\frac{T}{m}}\sigma(t_k, \bar X^m_{t_k})\,G_{k+1}, \; k=0,\ldots,m-1.
\end{equation}
Under the assumption
\begin{align}\exists \gamma  \!\in &(0,1],\;\exists C^{\rm Hol}_{b,\sigma},{\rm Lip}(b,\sigma)<\infty,\;\forall\, s,\, t \!\in[0,T], \; \forall\, x,\, y\!\in \R,\notag\\&|b(s,x)-b(t,y))| + |\sigma(s,x)-\sigma(t,y)| \le  C^{\rm Hol}_{b,\sigma}(1+|x|+|y|)|t-s|^{\gamma} + {\rm Lip}(b,\sigma)|x-y|,\label{eq:LipHol}
 \end{align}
 where, in addition to~\eqref{coeflip} for $b,\sigma$, a H\"older continuity of the coefficients in their time variable is supposed, the following standard error estimation holds for the scheme with non-truncated Brownian increments:
 \begin{equation}\label{errforteul}
 \Big \|\sup_{t\in [0,T]} |X_t-\bar X^m_t|  \Big\|_p\le \kappa^{(2)}_{b,\sigma,p}m^{-(\gamma\wedge \frac12)} (1+\|X_0\|_p).
\end{equation}
The next result addresses the convergence of { the version of the scheme} with truncated increments.

\begin{Proposition}\label{propconvforteul2}
Let us consider the truncated Euler scheme~\eqref{eq:EulerTrunc2} and let $p\in[1,+\infty)$. 

\smallskip
\noindent  {$(a)$} Under~\eqref{eq:LipHol}, there exists a finite constant $C_p$ such that for each initial random variable $X_0$ and each $m\ge 1$, the truncated Euler scheme~\eqref{eq:EulerTrunc} satisfies
   $$\Big \|\max_{k=0,\ldots,m}|X_{t_k}-\widetilde X^m_{t_k}|  \Big\|_p\le C_p\left(m^{-(\gamma\wedge \frac12)}+m^{\frac{1}{2p}}e^{-\frac{{\mathbf s}_m^2}{4p}}\right)(1+\|X_0\|_p).$$
  
  \smallskip
  \noindent  {$(b)$}  Moreover, if the coefficients $b,\sigma$ satisfy mere continuity on $[0,T]\times\R$ in addition to~\eqref{coeflip} and $X_0\in L^p(\P)$, then 
  $$
  \Big \|\max_{k=0,\ldots,m}|X_{t_k}-\widetilde X^m_{t_k}|  \Big\|_p\to 0\mbox{ as }m\to+\infty\;\mbox{  when }\;\liminf_m \frac{{\mathbf s}_m}{\sqrt{\ln m}}>\sqrt{2}.
  $$
\end{Proposition}

\noindent {\bf Proof of Proposition~\ref{propconvforteul2}.}
 {$(a)$} Under~\eqref{coeflip}, the standard estimation 
$$
\Big \|\sup_{t\in [0,T]} |X_t|  \Big\|_p+\Big \|\sup_{t\in [0,T]}|\bar X^m_t|  \Big\|_p +  \Big \|\max_{k=0,\ldots,m}|\widetilde X^m_{t_k}|  \Big\|_p\le \kappa^{(1)}_{b,\sigma,p,T}(1+\|X_0\|_p)
$$ 
holds. Reasoning like in the last step of the proof of Proposition~\ref{prop:ConvfortEuler}, one deduces that 
$$
\Big \|\max_{k=0,\ldots,m}|\bar X^m_{t_k} - \widetilde X^m_{t_k}|  \Big\|_p\le 2\kappa^{(1)}_{b,\sigma,2p,T} {\Big(\frac{m}{{\mathbf s}_m}e^{-\frac{{\mathbf s}_m^2}{2}}\Big)^{\frac{1}{2p}}}(1+\|X_0\|_p).
$$
Combining this inequality with~\eqref{errforteul} thanks to the triangle inequality, one obtains the first statement.

\smallskip
\noindent  {$(b)$ Since as $m\to +\infty$, when $\liminf_m\frac{{\mathbf s}_m}{\sqrt{\ln m}}>\sqrt{2}$, $\frac{m}{{\mathbf s}_m}e^{-\frac{{\mathbf s}_m^2}{2}}\to 0$}, to prove the second statement, it is enough to check that $\Big \|\sup_{t\in [0,T]} |X_t-\bar X^m_t|  \Big\|_p\to 0$. Let us suppose temporarily that $p\ge 2$ and $X_0\in L^{2p}(\P)$. Thanks to the spatial Lipschitz continuity of the coefficients, standard estimations based on the H\"older and Burkholder-Davis-Gundy inequalities and Gronwall's lemma ensure the existence of a finite constant not depending on $m$ and $X_0$ such that
\begin{align*}
  \Big \|\sup_{t\in [0,T]} |X_t-\bar X^m_t|  \Big\|^p_p&\le C\int_0^T\|b(s,X_s)-b({\underline s}, X_{\underline s})\|^p_p+\|\sigma(s,X_s)-\sigma({\underline s}, X_{\underline s})\|^p_pds.\end{align*}
With standard estimations on the moments of the increments of $X$ under~\eqref{coeflip}, we deduce that
\begin{align}\Big \|\sup_{t\in [0,T]} |X_t-\bar X^m_t|  \Big\|^p_p\le& C\int_0^T\|b(s,X_s)-b({\underline s}, X_{s})\|^p_p+\|\sigma(s,X_s)-\sigma({\underline s}, X_{s})\|^p_pds\notag\\&+C m^{-\frac p 2}(1+\|X_0\|_p)^p.\label{decomperreul}
\end{align}
Let now $K\in [1,+\infty)$ and $\varepsilon_{m,K}:=\sup_{0\le s\le t\le (s+\frac T m)\wedge T,x\in[-K,K]}(|b(s,x)-b(t,x)|+|\sigma(s,x)-\sigma(t,x)|)$. By uniform continuity of the coefficients $b,\sigma$ on the compact $[0,T]\times [-K,K]$, $\varepsilon_{m,K}\to 0$ as $m\to+\infty$. One has, using the affine growth of the coefficients in space uniform in time and $K\ge 1$,
\begin{align*}
   \sup_{s\in[0,T]}\left(|b(s,X_s)-b({\underline s}, X_{s})|+|\sigma(s,X_s)-\sigma({\underline s}, X_{s})|\right)&
  \le \varepsilon_{m,K}+C\sup_{t\in[0,T]}|X_t|{\mathbf 1}_{\{\sup_{t\in[0,T]}|X_t|\ge K\}}.
\end{align*}
Therefore, using the Cauchy-Schwarz inequality, then $\Big \|\sup_{t\in [0,T]} |X_t|  \Big\|_q\le \kappa^{(1)}_{b,\sigma,q}(1+\|X_0\|_q)$ for $q\in\{2p,1\}$ and the Markov inequality, we obtain
\begin{align}
  \sup_{s\in[0,T]}\big(\|b(s,X_s)-b({\underline s}, X_{s})\|^p_p&+\|\sigma(s,X_s)-\sigma({\underline s}, X_{s})\|^p_p\big)\notag\\&\le C\left(\varepsilon_{m,K}^p+\E\left[\sup_{t\in[0,T]}|X_t|^p{\mathbf 1}_{\{\sup_{t\in[0,T]}|X_t|\ge K\}}\right]\right)
  \notag\\&\le C\left(\varepsilon_{m,K}^p+  \Big\|\sup_{t\in[0,T]}|X_t|\Big\|_{2p}^p\P^{1/2}\left(\sup_{t\in[0,T]}|X_t|\ge K\right)\right)
  \notag\\&\le C\left(\varepsilon_{m,K}^p+  (1+\|X_0\|_{2p})^p\sqrt{1\wedge \frac{1+\|X_0\|_1}{K}}\right).\label{majodifftemps}
\end{align}

Plugging this estimation into~\eqref{decomperreul}, we obtain that when $p\ge 2$ and $X_0\in L^{2p}(\P)$,
\begin{align*}
   \Big \|\sup_{t\in [0,T]} |X_t-\bar X^m_t|  \Big\|^p_p\le C\left(m^{-\frac p 2}(1+\|X_0\|_p)^p+\varepsilon^p_{m,K}+  (1+\|X_0\|_{2p})^p\sqrt{1\wedge \frac{1+\|X_0\|_1}{K}}\right).
\end{align*}
Going to deterministic initial conditions $x_0\in\R$ and using that for $p\in[1,2]$, $\|\cdot\|^p_p\le\|\cdot\|^p_2$, we deduce with obvious notations, that when $p\ge 1$,
\begin{align*}
   \Big \|\sup_{t\in [0,T]} |X^{x_0}_t-\bar X^{m,x_0}_t|  \Big\|^p_p\le C\left(m^{-\frac p 2}(1+|x_0|)^p+\varepsilon^p_{m,K}+  (1+|x_0|)^p\left(1\wedge \frac{1+|x_0|}{K}\right)^{\frac{1}{2}\wedge\frac{p}{4}}\right).
\end{align*}
By Blagove$\check{\rm   s}\check{\rm  c}$enkii-Freidlin's theorem~\cite[Theorem 13.1, Section V.12-13, p.136]{RogersWilliamsII}, there exists a measurable function $F:\R\times C([0,T],\R)\to C([0,T],\R)$ continuous in its first variable such that $X^{x_0}=F(x_0,W)$ and $X=F(X_0,W)$. Moreover, by the inductive definition of the Euler scheme, there also exists a measurable function $F_m: \R\times C([0,T],\R)\to C([0,T],\R)$ such that $\bar X^{m,x_0}=F_m(x_0,W)$ and $\bar X^{m}=F_m(X_0,W)$. Therefore integrating the last inequality in $x_0$ with respect to the law of $X_0$ and using the independence of $X_0$ and $W$, we deduce that when $p\ge 1$ and $X_0\in L^{p}(\P)$
\begin{align*}
   \Big \|\sup_{t\in [0,T]} |X_t-\bar X^m_t|  \Big\|_p\le C\left(m^{-\frac1 2}(1+\|X_0\|_p)+\varepsilon_{m,K}+\E^{1/p}\left[(1+|X_0|)^p\left(1\wedge \frac{1+|X_0|}{K}\right)^{\frac{1}{2}\wedge\frac{p}{4}}\right]\right).
\end{align*}
Since, when $X_0\in L^p(\P)$,  the expectation in the right-hand side goes to $0$ as $K\to+\infty$ by Lebesgue's theorem, we can choose $K$ so that the last term in the right-hand side is arbitrarily small and then the sum of the two first terms goes to $0$ when $m\to+\infty$. The conclusion follows.
\hfill$\Box$

\subsection{Ordering of the Euler schemes with truncated noise}

In the next proposition, we focus on  orderings of the truncated Euler schemes. To deal at the same time with the schemes~\eqref{eq:EulerTrunc} and~\eqref{eq:EulerTrunc2} we set $(\widetilde X^m_0,\widetilde X^m_0)=(X_0,Y_0)$ and
\begin{align}
  &\widetilde X^m_{t_{k+1}}= \widetilde  X^m_{t_k} + \frac T m\widetilde b^m_k(\widetilde  X^m_{t_k}) + \sqrt{\frac Tm}\widetilde \sigma^m_k(\widetilde   X^m_{t_k})\,Z^m_{k+1},\notag\\
  &\widetilde Y^m_{t_{k+1}}= \widetilde  Y^m_{t_k} + \frac Tm\widetilde \beta^m_k(\widetilde  Y^m_{t_k})ds + \sqrt{\frac Tm}\widetilde \vartheta^m_k(\widetilde   Y^m_{t_k})\,Z^m_{k+1},\; k=0,\ldots,m-1.\label{2eulertrunc}
\end{align}
where  either

\smallskip
 for $k=0,\ldots,m-1$, $(\widetilde b^m_k(\cdot),\widetilde\sigma^m_k(\cdot),\widetilde\beta^m_k(\cdot),\widetilde\varsigma_k^m(\cdot))$ is equal to $$\left(\frac{m}{T}\int_{t_k}^{t_{k+1}}b(s,\cdot)ds,\sqrt{\frac{m}{T}\int_{t_k}^{t_{k+1}}\sigma^2(s,\cdot)ds},\frac{m}{T}\int_{t_k}^{t_{k+1}}\beta(s,\cdot)ds,\sqrt{\frac{m}{T}\int_{t_k}^{t_{k+1}}\vartheta^2(s,\cdot)ds}\right)$$ (schemes~\eqref{eq:EulerTrunc})  
 
 \smallskip
 or
 
 \smallskip
  for $k=0,\ldots,m-1$, $(\widetilde b^m_k(\cdot),\widetilde\sigma^m_k(\cdot),\widetilde\beta^m_k(\cdot),\widetilde\varsigma_k^m(\cdot))=\left(b(t_k,\cdot),\sigma(t_k,\cdot),\beta(t_k,\cdot),\vartheta(t_k,\cdot)\right)$ (schemes~\eqref{eq:EulerTrunc2}).
For $x,y\in\R$, we also denote by $(\widetilde X^{m,x}_{t_k})_{k=0,\hdots,m}$ and $(\widetilde Y^{m,y}_{t_k})_{k=0,\hdots,m}$ these truncated Euler schemes started from $\widetilde X^m_{0}=x$ and $\widetilde Y^m_{0}=y$.

Let \begin{align}\label{eq:conv_sigma}
&\mathbf{a}_\sigma = \inf\big\{a\ge 0 : x\mapsto \sigma^2(t,x) +ax^2 \mbox{ is convex for all } \; t\!\in [0,T]\},\\\label{eq:nondec_b}
      &\mathbf{c}_b = \inf\big\{c\ge 0 : x\mapsto b(t,x) +cx \mbox{ is non-decreasing for all } \; t\!\in [0,T]\},
\end{align}
and {let} $\mathbf{a}_{\vartheta}$, $\mathbf{c}_\beta$ be defined likewise. Let also $\mathbf{c}_\sigma:=\mathbf{a}_\sigma+ ({\rm Lip}(\sigma))^2$, 
$\mathbf{m}_{b,\sigma}:=\left(\frac{2\mathbf{c}_b}{\sqrt{\mathbf{c}_\sigma+2\mathbf{c}_b}-\sqrt{\mathbf{c}_\sigma}}\right)^2T$ and {let}  $\mathbf{c}_\vartheta,\mathbf{m}_{\beta,\vartheta}$ be defined likewise.
\begin{Remark}
  The finiteness of $\inf\big\{a\ge 0 : x\mapsto f(x) +ax^2\big\}$ is often known as the {\em semi-convexity} property of $f:\R\to\R$. It has been seemingly introduced in Mathematical Finance in~\cite{Caverhilletal1990} and can also be found in various papers on the pricing and hedging of American style options like~\cite{Lamb2002},~\cite{BalPag2003} or~\cite{Ballyetal2005}. {A typical example of semi-convex function is a differentiable function with Lipschitz differential.} 
\end{Remark}

\begin{Proposition}\label{prop:DirConvforEuler} Assume that the functions $b$, $\beta$, $\sigma$, $\vartheta : \R\times [0,T]\to \R$ satisfy~\eqref{coeflip}. Also assume 
$$
(\mathcal X) \equiv \left\{\begin{array}{ll} (i) & \mathbf{a}_{\sigma}<+\infty  \\
(ii)& \forall\, t\!\in [0,T],\; b(t,\cdot) \mbox{ convex}
\end{array}\right.  \quad \mbox{ or }\quad (\mathcal Y) \equiv\left\{\begin{array}{ll} (i) & \mathbf{a}_{\vartheta} <+\infty \\
(ii) &\forall\, t\!\in [0,T],\;  \beta(t,\cdot) \mbox{ convex}.
\end{array}\right. 
$$
We consider the Euler schemes $\widetilde X^m$ and $\widetilde Y^m$ given in~\eqref{2eulertrunc} starting respectively from $X_0$ and $Y_0\!\in L^1(\P)$, with $m\in \N$, $m\ge {\mathbf m}_{b,\sigma}$ and truncated noise at level 
$\mathbf{s}_m\in\left[0,\frac{\sqrt{m}}{2{\rm Lip}(\sigma)\sqrt{T}}\right]$ if $(\mathcal X)$ holds or with $m\in \N$, $m\ge {\mathbf m}_{b,\sigma}$ and truncated noise at level $\mathbf{s}_m\in\left[0,\frac{\sqrt{m}}{2{\rm Lip}(\vartheta)\sqrt{T}}\right]$ if $(\mathcal Y)$ holds.

\smallskip
\noindent $(a)$ {\em Non-decreasing directional convex ordering}.  We have 
\begin{align*}
   \Big( b \le \beta \mbox{ and }0\le \sigma\le \vartheta &\mbox{ on }[0,T]\times \R\mbox{ and }X_0\preceq_{icv} Y_0\Big)\Longrightarrow   (\widetilde X^{m}_{t_k})_{k=0,\ldots,m}\le_{diricv}(\widetilde Y^{m}_{t_k})_{k=0,\ldots,m}.
\end{align*}
Moreover, when  $(\mathcal X)$ (resp. $(\mathcal Y)$) holds, then for each directionally convex function $f:\R^{m+1}\to \R$ non-decreasing in each of its variables,
$$
x\mapsto \E \,f((\widetilde X^{m,x}_{t_k})_{k=0,\ldots,m})\in\R\cup\{+\infty\} \left(\mbox{resp. }y\mapsto \E \,f((\widetilde Y^{m,y}_{t_k})_{k=0,\ldots,m})\in\R\cup\{+\infty\}\right)$$ is non-decreasing and convex.

\smallskip
\noindent $(b)$ {\em Directional convex ordering}. Let  $b(t,x)=\beta(t,x) = \lambda(t)+\mu(t)x$, $(t,x)\!\in [0,T]\times \R$  be a time dependent affine function in space where $\lambda$ and $\mu$ are bounded Borel functions. 
Then
\begin{align*}
   \Big(0\le \sigma\le \vartheta&\mbox{ on } [0,T]\times \R \mbox{ and } X_0\preceq_{cvx} Y_0\Big)\Longrightarrow (\widetilde X^{m}_{t_k})_{k=0,\ldots,m}\le_{dircvx}(\widetilde Y^{m}_{t_k})_{k=0,\ldots,m}.
\end{align*}
Moreover, when  $(\mathcal X)$ (resp. $(\mathcal Y)$) holds, then for each directionally convex function $f:\R^{m+1}\to \R$,
$$
x\mapsto \E \,f((\widetilde X^{m,x}_{t_k})_{k=0,\ldots,m})\in\R\cup\{+\infty\} \left(\mbox{resp. }y\mapsto \E \,f((\widetilde Y^{m,y}_{t_k})_{k=0,\ldots,m})\in\R\cup\{+\infty\}\right)$$ is convex.
\end{Proposition}

\begin{Remark}\label{remdervar} In view of Property {\bf P3}  in Remark~\ref{remdircon}, $(b)$ (resp. $(a)$) implies that when  $(\mathcal X)$ holds and $b$ is affine (resp. convex) in space, then $x\mapsto \E g(\widetilde X^{m,x}_{T})$ is convex (resp. non-decreasing and convex) when $g:\R\to\R$ is convex (resp. non-decreasing and convex).
\end{Remark}

\medskip The proof of Proposition~\ref{prop:DirConvforEuler} relies on the  next results which deals with directional convexity propagation by a truncated Euler operator.

\begin{Proposition}\label{prop:propconv}

Assume that 
  
\smallskip
  \begin{itemize}
  \item the function $\beta:\R\to\R$ is convex and such that 
\begin{equation}\label{eq:beta}
 c_\beta:=\inf\big\{c\ge 0:\;\R\ni x\mapsto \beta(x)+cx\mbox{ is non-decreasing}\big\}<+\infty.
\end{equation}
  \item the function $\varsigma:\R\to\R_+$ is Lipschitz continuous with constant ${\rm Lip}(\varsigma)$ and such that
\begin{equation}\label{eq:a_sigma}
  a_{\varsigma} := \inf \big\{a\ge 0:x\mapsto \varsigma^2(x) +a x^2 \hbox{ is convex}\big\}<+\infty,
\end{equation}
and let $c_\varsigma:=a_\varsigma+ ({\rm Lip}(\varsigma))^2$.
\end{itemize} 
For $h\in\R_+^*$, let $\mathcal{E}^h :\R\times \R  \to\R$ be the function defined by  
$$
\mathcal{E}^h(x,z):= x+\sqrt{h}\,\varsigma(x)z+h\beta(x),
$$ 
and $Z^s={\mathbf 1}_{\{G\in[-s,s]\}}G$ with $G\sim{\cal N}_1(0,1)$ and $s\in \left[0,\frac{1}{2{\rm Lip}(\varsigma)\sqrt{h}}\right]$. Then, the following holds: 

\medskip
\noindent  $(i)$ when $h\le \frac1{2c_\beta}$, the random function $\mathcal{E}^h(\cdot, Z^s)$ is non-decreasing when the codomain is endowed with the stochastic order,

\smallskip
\noindent $(ii)$ Let $h\le \bar{h}_{\beta,\varsigma}:=\left(\frac{\sqrt{c_\varsigma+2c_\beta}-\sqrt{c_\varsigma}}{2c_\beta}\right)^2$, $k\in\N$ and $f:\R^{k+2}\to\R$.

If $f$ is Lipschitz continuous, directionally convex and non-decreasing in each of its variables, then $\R^k\times\R\ni (y,x)\mapsto\E \,f\big(y,x,\mathcal{E}^h(x, Z^s)\big)$ is Lipschitz continuous, directionally convex and non-decreasing in each of its variables.

If $\beta$ is affine and $f$ is Lipschitz continuous and directionally convex, then $\R^k\times\R\ni (y,x)\mapsto\E \,f\big(y,x,\mathcal{E}^h(x, Z^s)\big)$ is Lipschitz continuous and directionally convex.

\smallskip
We use the conventions $\frac{1}{0}=+\infty$ in particular in the definition of the interval where the threshold $s$ can take its values and  $\frac{\sqrt{c_\varsigma+2c_\beta}-\sqrt{c_\varsigma}}{2c_\beta}$ {is conventionally equal to} $\lim_{x\to 0+}\frac{\sqrt{c_\varsigma+x}-\sqrt{c_\varsigma}}{x}=\frac{1}{2\sqrt{c_\varsigma}}$ when $c_\beta=0$.
\end{Proposition}
\begin{Remark} 
When $Z$ is centered and such that ${\rm Var}(Z)\in (0,+\infty)$, a necessary condition for $\mathcal{E}^h(\cdot, Z)$ to be convex {when the codomain is endowed with the convex order} is
 \begin{align*}
   \forall x,y\in\R,\;\forall \alpha\in[0,1],\;&\E\,\mathcal{E}^h(\alpha x+(1-\alpha)y, Z)=\alpha\E\,\mathcal{E}^h(x, Z)+(1-\alpha)\E\,\mathcal{E}^h(y, Z)\\&\mbox{ and }\E\,[\mathcal{E}^h(\alpha x+(1-\alpha)y, Z)^2]\le \alpha\E\,[\mathcal{E}^h(x, Z)^2]+(1-\alpha)\E\,[\mathcal{E}^h(y, Z)^2]
 \end{align*}
 which also writes $\beta(\alpha x+(1-\alpha)y)=\alpha\beta(x)+(1-\alpha)\beta(y)$ and 
 $$\varsigma^2(\alpha x+(1-\alpha)y)+\frac{(\alpha x+(1-\alpha)y)^2}{h{\rm Var}(Z)}\le \alpha\left(\varsigma^2(x)+\frac{x^2}{h{\rm Var}(Z)}\right)+(1-\alpha)\left(\varsigma^2(y)+\frac{y^2}{h{\rm Var}(Z)}\right)$$
 so that $\beta$ is affine and $a_\varsigma\le \frac{1}{h{\rm Var}(Z)}$.
\end{Remark}
\begin{Remark}\label{rem:contrex}
   Let $Q(x,dy)$ denote the distribution of $x+\sqrt{h}\varsigma(x)Z^s$ with $s\in\left[0,\frac{1}{2{\rm Lip}(\varsigma)\sqrt{h}}\right]$. When $h\le \frac{1}{4c_\varsigma}$, the Markov kernel $Q$ is such that when $\varphi:\R^2\to\R$ is directionally convex, then $x\mapsto \E[\varphi(x,x+\sqrt{h}\varsigma(x)Z^s)]=\int_\R \varphi(x,y)Q(x,dy)$ is convex by Proposition \ref{prop:propconv} applied with $k=0$ and $\beta\equiv 0$. But this convexity property may fail when $\varphi$ is convex . 
 
 Indeed, restricted to the class of functions $\varphi(x,y)=f(y-x)$ with  convex $f:\R\to\R$ for which $\E[\varphi(x,x+\sqrt{h}\varsigma(x)Z^s)]=\E[f(\sqrt{h}\varsigma(x)Z^s)]$, the convexity preservation would imply that 
 $$
 x\mapsto  \,{\E\,f(\sqrt{h}\varsigma(x)Z^s)}\quad \mbox{  is convex.}
 $$
  Setting $f(x)= |x|$ (which is also the choice made in the proof of Proposition \ref{prop:CNconv2}) implies that 
   $$
 x\mapsto  \sqrt{h}\,{|\varsigma(x)|\, \E\,|Z^s|}\quad \mbox{  is convex.}
 $$
When $\E\, |Z^s|>0$, this implies that $\varsigma = |\varsigma|$ is convex. 
 
\smallskip
As a consequence, the proof of~\cite[Proposition 3.1]{Ruschendorf2}, restated in~\cite[Lemma 1]{BRADAC} and in~\cite[p54-55]{Ruschendorf3}, is not valid when ${\cal F}$ is the class of convex functions. 
\end{Remark}

\noindent {\bf Proof of Proposition~\ref{prop:DirConvforEuler}.}  For $k=0,\ldots,m-1$, and $f_{k+1}:\R^{k+2}\to\R$ Lipschitz continuous, we set
\begin{align*}
  \R^{k+1}\ni x_{0:k}=(x_0,\ldots,x_k)\mapsto \widetilde P_k^{(b,\sigma),m}f_{k+1}(x_{0:k})&=\E \,f_{k+1}\left(x_{0:k},x_k+\frac T m\widetilde b^m_k(x_k) + \sqrt{\frac Tm}\widetilde \sigma^m_k(x_k)\,Z^m_{k+1}\right)\\
  \R^{k+1}\ni x_{0:k}\mapsto \widetilde P_k^{(\beta,\vartheta),m}f_{k+1}(x_{0:k})&=\E \,f_{k+1}\left(x_{0:k},x_k+\frac T m\widetilde \beta^m_k(x_k) + \sqrt{\frac Tm}\widetilde \vartheta^m_k(x_k)\,Z^m_{k+1}\right).
\end{align*}
Since $b,\sigma,\beta,\vartheta$ satisfy~\eqref{coeflip}, $\widetilde P_k^{(b,\sigma),m}f_{k+1}$ and $\widetilde P_k^{(\beta,\vartheta),m}f_{k+1}$ are Lipschitz continuous.

Still for $k=0,\ldots,m-1$, we set $\widetilde P_{k:m-1}^{(b,\sigma),m}=\widetilde P^{(b,\sigma),m}_k \circ \cdots\circ \widetilde P_{m-1}^{(b,\sigma),m}$ and $\widetilde P_{k:m-1}^{(\beta,\vartheta),m}=\widetilde P^{(\beta,\vartheta),m}_k \circ \cdots\circ \widetilde P_{m-1}^{(\beta,\vartheta),m}$. By the Markov property, for $f:\R^{m+1}\to\R$ Lipschitz continuous, we have 
\begin{align*}
   \forall x\in\R,\;&\E\, f\big((\widetilde X^{m,x}_{t_k})_{k=0,\ldots,m}\big) = \widetilde P^{(b,\sigma),m}_{0:m-1}f(x)\mbox{ and }\E\, f\big((\widetilde Y^{m,x}_{t_k})_{k=0,\ldots,m}\big) = \widetilde P^{(\beta,\vartheta),m}_{0:m-1}f(x),\\
&\E\, f\big((\widetilde X^{m}_{t_k})_{k=0,\ldots,m}\big)=\E\widetilde P^{(b,\sigma),m}_{0:m-1}f(X_0)\mbox{ and }\E\, f\big((\widetilde Y^{m}_{t_k})_{k=0,\ldots,m}\big) = \E\widetilde P^{(\beta,\vartheta),m}_{0:m-1}f(Y_0).
\end{align*}
Let $f:\R^{m+1}\to\R$ be Lipschitz continuous, non-decreasing in each of its variables and directionally convex. We assume $(\mathcal X)$, that  $m\ge {\mathbf m}_{b,\sigma}$ and
$\mathbf{s}_m\in\left[0,\frac{\sqrt{m}}{2{\rm Lip}(\sigma)\sqrt{T}}\right]$. Then for each $k=0,\ldots,m-1$, $a_{\widetilde\sigma^m_k}\le \mathbf{a}_\sigma$, ${\rm Lip}(\widetilde\sigma^m_k)\le {\rm Lip}(\sigma)$ so that $c_{\widetilde\sigma^m_k}\le \mathbf{c}_\sigma$ and  $c_{\widetilde b^m_k}\le \mathbf{c}_b$. Since $\R^*_+\times \R^*_+\ni(z,w)\mapsto \frac{\sqrt{z+w}-\sqrt{z}}{w}$ is non-increasing in each variable, $\left(\frac{\sqrt{c_{\widetilde\sigma^m_k}+2c_{\widetilde b^m_k}}-\sqrt{c_{\widetilde\sigma^m_k}}}{2c_{\widetilde b^m_k}}\right)^2\ge \left(\frac{\sqrt{\mathbf{c}_\sigma+2\mathbf{c}_b}-\sqrt{\mathbf{c}_\sigma}}{2\mathbf{c}_b}\right)^2$ and $h=\frac{T}{m}\le \bar h_{\widetilde b^m_k,\widetilde \sigma^m_k}$. Combining Proposition~\ref{prop:propconv} applied with $(\beta,\varsigma)=(\widetilde b^m_k,\widetilde\sigma^m_k)$, $h=\frac T m$, $s=\mathbf{s}_m\in\left[0,\frac{1}{2{\rm Lip}(\widetilde\sigma^m_k)\sqrt{h}}\right]$ with an obvious backward inductive reasoning, we check that for all $k=0,\ldots,m-1$, $\widetilde P_{k:m-1}^{(b,\sigma),m}f$ is Lipschiz continuous, non-decreasing  in each of its variables and directionally convex. In particular, $\R\ni x\mapsto \E\, f\big((\widetilde X^{m,x}_{t_k})_{k=0,\ldots,m}\big) = \widetilde P^{(b,\sigma),m}_{0:m-1}f(x)$ is non-decreasing and convex.
We now proceed by a backward induction on $k$ to show that $\widetilde P ^{(b,\sigma),m}_{k:m-1} f \le \widetilde P_{k:m-1}^{(\beta,\vartheta),m}f$ for $k=0\ldots,m$, where, by convention, $\widetilde  P^{(b,\sigma),m}_{m:m-1}f=\widetilde  P^{(\beta,\vartheta),m}_{m:m-1}f=f$ so that the inequality holds for $k=m$.
Let for $f_{k+1}:\R^{k+2}\to\R$ Lipschitz continuous and $(x_{0:k},u,v)\in\R^{k+1}\times\R\times\R$, $Q_kf_{k+1}(x_{0:k},u,v)= \E\, f_{k+1}\left(x_{0:k},x_k+\frac T m u +v\sqrt{\frac T m}Z^m_{k+1}\right)$. 

 Assume that $\widetilde P ^{(b,\sigma),m}_{k+1:m-1} f \le \widetilde P_{k+1:m-1}^{(\beta,\vartheta),m}f$ for some $k\in\{0,\ldots,m-1\}$. Then, since $\widetilde P^{(b,\sigma)}_{k+1:m-1} f$ is non-decreasing in its last variable, $\forall (x_{0:k},v)\in\R^{k+1}\times \R$, $\R\ni u\mapsto Q_k\widetilde P^{(b,\sigma)}_{k+1:m-1} f(x_{0:k},u,v)$ is non-decreasing. Moreover, since $\widetilde P^{(b,\sigma)}_{k+1:m-1} f$ is convex in its last variable, $\forall (x_{0:k},u)\in\R^{k+1}\times \R$, $\R\ni v\mapsto Q_k\widetilde P^{(b,\sigma)}_{k+1:m-1} f(x_{0:k},u,v)$ is convex and, by Jensen's inequality, minimal for $v=0$, so that $\R_+\ni v\mapsto Q_k\widetilde P^{(b,\sigma)}_{k+1:m-1} f(x_{0:k},u,v)$ is non-decreasing.
Combined with $\widetilde b^m_k(x_k)\le\widetilde\beta^m_k(x_k)$ and $0\le \widetilde\sigma^m_k(x_k)\le\widetilde\vartheta^m_k(x_k)$, this ensures that for $x_{0:k}\in\R^{k+1}$,
\begin{align*}
\widetilde  P^{(b,\sigma),m}_{k:m-1} f(x_{0:k}) =  Q_k\widetilde  P^{(b,\sigma),m}_{k+1:m-1} f\big( x_{0:k},\widetilde b^m_k(x_k),\widetilde \sigma^m_k(x_k)\big)
  &\le Q_k\widetilde  P^{(b,\sigma),m}_{k+1:m-1} f \big( x_{0:k},\widetilde \beta^m_k(x_k),\widetilde\vartheta^m_k(x_k)\big).\end{align*}
With the positivity of $Q_k$ as a linear operator on functions and the induction hypothesis, we deduce that
$$\widetilde  P^{(b,\sigma),m}_{k:m-1} f(x_{0:k})\le Q_k\widetilde   P_{k+1:m-1}^{(\beta,\vartheta),m} f\big( x_{0:k},\widetilde \beta^m_k(x_k),\widetilde\vartheta^m_k(x_k)\big) = \widetilde P^{(\beta,\vartheta),m}_{k:m-1} f(x_{0:k}),$$
which concludes the proof by induction.
The monotonicity and convexity of  $\widetilde P^{(b,\sigma),m}_{0:m-1}f$ combined with $X_0\preceq_{icv} Y_0$ and the previous inequality for $k=0$ then imply that
$$ 
\E\, f\big((\widetilde X^{m}_{t_k})_{k=0,\ldots,m}\big)=\E\,  \widetilde P^{(b,\sigma),m}_{0:m-1}f(X_0)\le \E\,  \widetilde P^{(b,\sigma),m}_{0:m-1}f(Y_0)\le \E\,  \widetilde P^{(\beta,\vartheta),m}_{0:m-1}f(Y_0) = \E\, f\big((\widetilde Y^{m}_{t_k})_{k=0,\ldots,m}\big). 
$$
When $\big({\cal Y}\big)$ holds and $m\ge {\mathbf m}_{\beta,\vartheta}$, $\mathbf{s}_m\in\left[0,\frac{\sqrt{m}}{2{\rm Lip}(\vartheta)\sqrt{T}}\right]$, then for each $k=0,\ldots,m-1$, $x_{0:k}\mapsto P^{(\beta,\vartheta),m}_{k:m-1}f(x_{0:k})$ is non-decreasing and convex by Proposition~\ref{prop:propconv} and backward induction and we use
\begin{align*}
   Q_k\widetilde  P^{(b,\sigma),m}_{k+1:m-1} f\big( x_{0:k},\widetilde b^m_k(x_k),\widetilde \sigma^m_k(x_k)\big)&\le Q_k\widetilde  P^{(\beta,\vartheta),m}_{k+1:m-1} f \big( x_{0:k},\widetilde b^m_k(x_k),\widetilde \sigma^m_k(x_k)\big)\\&\le Q_k\widetilde   P_{k+1:m-1}^{(\beta,\vartheta),m} f\big( x_{0:k},\widetilde \beta^m_k(x_k),\widetilde\vartheta^m_k(x_k)\big)
\end{align*}
in the backward inductive proof of $\widetilde P ^{(b,\sigma),m}_{k:m-1} f \le \widetilde P_{k:m-1}^{(\beta,\vartheta),m}f$ for $k=0,\ldots,m-1$. We also use $\E\,  \widetilde P^{(b,\sigma),m}_{0:m-1}f(X_0)\le \E\,  \widetilde P^{(\beta,\vartheta),m}_{0:m-1}f(X_0)\le \E\,  \widetilde P^{(\beta,\vartheta),m}_{0:m-1}f(Y_0)$ to conclude that 
$$
\E\, f\big((\widetilde X^{m}_{t_k})_{k=0,\ldots,m}\big)\le \E\, f\big((\widetilde Y^{m}_{t_k})_{k=0,\ldots,m}\big).
$$
 With Lemma~\ref{lem:approxdirconv}, we deduce that $(\widetilde X^{m}_{t_k})_{k=0,\ldots,m}\le_{diricv}(\widetilde Y^{m}_{t_k})_{k=0,\ldots,m}$.

Let us now check that mere directional convexity and monotonicity in each variable of $f$ (without Lipschitz continuity) is enough to ensure that $x\mapsto \E\, f\big((\widetilde X^{m,x}_{t_k})_{k=0,\ldots,m}\big)$ is non-decreasing and convex. By Lemma~\ref{lem:approxdirconv},  there exists a non-decreasing sequence of Lipschitz continuous, non-decreasing in each variable and directionally convex functions  $(f_n)_{n\ge 1}$ converging pointwise to $f$. By Beppo-Levi's theorem, $ \E\,\big[(f_n-f_1)\big((\widetilde X^{m,x}_{t_k})_{k=0,\ldots,m}\big)\big]$ converges to $\E\,\big[(f-f_1)\big((\widetilde X^{m,x}_{t_k})_{k=0,\ldots,m}\big)\big]$ possibly equal to $+\infty$ as $n\to+\infty$. Since $f_1$ is Lipschitz continuous, $\E\,\big|f_1\big((\widetilde X^{m,x}_{t_k})_{k=0,\ldots,m}\big)\big|<\infty$ and we deduce that $x\mapsto \E\,f_n\big((\widetilde X^{m,x}_{t_k})_{k=0,\ldots,m}\big)$ converges pointwise to $x\mapsto \E\,f\big((\widetilde X^{m,x}_{t_k})_{k=0,\ldots,m}\big)$ as $n\to+\infty$ so that the non-decreasing convex property of the former functions is transferred to the latter.

The proof is analogous for the convex ordering using that the affine property of the drift coefficients enables backward propagation of directional convexity and their equality permits to freeze $u=\widetilde b^m_k(x_k)=\widetilde \beta^m_k(x_k)$ in $Q_kf_{k+1}(x_{0:k},u,v)$.\hfill$\Box$

\bigskip
\noindent {\bf Proof of Proposition~\ref{prop:propconv}.}   
\noindent   {\sc Step~1}. ({\em Preliminaries}). 

Let $\rho$ be a $C^\infty$  probability density on the real line with compact support such that \begin{equation}
   \int_{\R} u\r(u)du = 0\mbox{ and }\int u^2\rho(u)du =1.\label{eqcentred}
\end{equation} We associate to $\rho$ its sequence of mollifiers  $\rho_n(x)=n\rho(nx)$, $n\in \N$.

We set, for every $n\ge 1$,  
$$
\varsigma_n(x)=\sqrt{\frac1n+\rho_n\star \varsigma^2(x)}\quad\mbox{ and }\quad \beta_n=\rho_n\star \beta.
$$
The convex function $\beta$ being continuous, the function $\beta_n$ is infinitely differentiable, convex and such that $x\mapsto \beta_n(x)+c_\beta x$ is non-decreasing, so that
\begin{equation}
   \beta'_n\ge -c_\beta.\label{minoderbetn}
\end{equation} 
Furthermore, $\beta_n$ converges pointwise to $\beta$ as $n\to+\infty$. By continuity of {$\varsigma^2$,  $\rho_n\star\varsigma^2$ converges to $\varsigma^2$
so that} $\varsigma_n$ also converges pointwise to $\varsigma$. As a consequence,
$$
\mathcal{E}^h_n(x,Z^s):=x+\sqrt{h}\,\varsigma_n(x)Z^s+h\beta_n(x)\to\mathcal{E}^h(x,Z^s)\quad\mbox{ as }\quad n\to+\infty.
$$
Moreover $$\big|\mathcal{E}^h_n(x,Z^s)-\mathcal{E}^h(x,Z^s)\big|\le \sqrt{h}|\varsigma_n(x)-\varsigma(x)||Z|+h|\beta_n(x)-\beta(x)|\to 0\mbox{ as }n\to+\infty.$$
For the non-decreasing stochastic order, this ensures that the non-decreasing property of $x\mapsto\mathcal{E}^h_n(x,Z^s)$ that we are going to check for each $n\in\N$ is transferred to $x\mapsto\mathcal{E}^h(x,Z^s)$. This also ensures that for a Lipschitz function $f:\R^{k+2}\to\R$ and $(y,x)\in\R^k\times\R$, one has $\E \,f\big(y,x,\mathcal{E}^h_n(x,Z^s)\big)\to\E \,f\big(y,x,\mathcal{E}^h(x,Z^s)\big)$ as $n\to+\infty$. Therefore it is enough to prove that the conclusion holds for $\E \,f\big(y,x,\mathcal{E}^h_n(x,Z^s)\big)$ for each $n\in\N$, which we are also going to do.

Since, by Jensen's inequality $\varsigma \star \rho_n\le \sqrt{\varsigma^2\star \rho_n}\le \varsigma_n$, we obtain that for $x,x'\in \R$,
\begin{align*}
   |\varsigma^2_n(x)-\varsigma^2_n(x')|&\le\int_\R|\varsigma(x-z)-\varsigma(x'-z)|(\varsigma(x-z)+\varsigma(x'-z))\rho_n(z)dz\\&\le {\rm Lip}(\varsigma)|x-x'|\left(\varsigma \star \rho_n(x)+\varsigma \star \rho_n(x')\right)\le {\rm Lip}(\varsigma)|x-x'|(\varsigma_n(x)+\varsigma_n(x')).
\end{align*}
Hence the function $\varsigma_n$ is Lipschitz continuous with constant ${\rm Lip}(\varsigma)$. Moreover, since the Lipschitz function $\varsigma$ has at most affine growth, the non-negative function $\rho_n\star\varsigma^2$ is infinitely differentiable and so is $\varsigma_n$. The Lipschitz property of $\varsigma_n$ implies that 
\begin{equation}
   |\varsigma_n'|\le {\rm Lip}(\varsigma).\label{lipvarsig}
\end{equation}

 Moreover,   it is  clear that $\varsigma (x)^2 +a_{\varsigma}x^2$ is convex since the infimum holds as a minimum. 
Hence its convolution by $\rho_n$ is convex too and one checks that it is infinitely differentiable. On the other hand, using~\eqref{eqcentred},  one has
\[
\int\big[\varsigma^2 (x-y) +a_{\varsigma}(x-y)^2\big]\rho_n(y)dy = \rho_n\star \varsigma^2(x) +a_{\varsigma} x^2 + \tfrac{a_{\varsigma}}{n^2}=\varsigma^2_n(x)+a_{\varsigma} x^2 + \tfrac{a_{\varsigma}}{n^2}-\tfrac 1n,
\]
so that $\varsigma^2_n (x)+ a_{\varsigma}x^2$ is convex and  $(\varsigma_n^2)''\ge -2\,a_{\varsigma}$. As a consequence, \begin{align} \forall x\in\R,\;\varsigma_n\varsigma_n''(x)&=\frac12(\varsigma_n^2)''(x)-(\varsigma'_n(x))^2\ge -a_\varsigma-({\rm Lip}(\varsigma))^2=-c_\varsigma.\label{minosigsig''}
\end{align}

\smallskip
\noindent {\sc Step 2} ({\em Non-decreasing stochastic ordering}). 
 When $h\le \frac{1}{2c_\beta}$, then $h\beta_n'(x)\ge -\frac12$ and, with $s\in\left[0,\frac{1}{2{\rm Lip}(\varsigma)\sqrt{h}}\right]$ and~\eqref{lipvarsig}, we deduce that
 \begin{equation}
   \forall x\in\R,\;\partial_x\mathcal{E}^h_n(x,Z^s)=1+\sqrt{h}\varsigma_n'(x)Z^s+h\beta_n'(x)\ge 1- \frac12-\frac12  =0.\label{derxcale}
 \end{equation}
Therefore, $x\mapsto \mathcal{E}^h_n(x,Z^s)$ is non-decreasing for the non-decreasing stochastic ordering and $\R^k\times\R\ni (y,x)\mapsto \E \,f\big(y,x,\mathcal{E}^h_n(x,Z^s)\big)$ is non-decreasing in each of its variables when $f:\R^{k+2}\to\R$ is non-decreasing in each of its variables and Lipschitz continuous.

\medskip
\noindent {\sc Step 3} ({\em Differentiation}). Let $f:\R^{k+2}\to\R$ be continuously differentiable with bounded first order derivatives and $\R^k\times \R\ni(y,x)\mapsto g(y,x)=\E \,f\big(y,x,\mathcal{E}^h_n(x,Z^s)\big)$. One has 
\begin{align}
   \partial_x &g(y,x)=\E\,\partial_{k+1}f(y,x,\mathcal{E}^h_n(x,Z^s))+\E\left[\partial_{k+2}f(y,x,\mathcal{E}^h_n(x,Z^s))(1+\sqrt{h}\varsigma_n'(x)Z^s+h\beta_n'(x))\right].\label{derxg}
\end{align}
And when $f$ is twice continuously differentiable with bounded first and second order derivatives, then 
for $i,j\in\{1,\cdots,k\}$,
\begin{align}
&  \partial_{y_iy_j}g(y,x)=\E\,\partial_{ij}f\big(y,x,\mathcal{E}^h_n(x,Z^s)\big),\label{dersecyy}\\
& \partial_{y_ix}g(y,x)=\E\,\partial_{ik+1}f\big(y,x,\mathcal{E}^h_n(x,Z^s)\big)+\E\left[\partial_{ik+2}f(y,x,\mathcal{E}^h_n(x,Z^s))(1+\sqrt{h}\varsigma_n'(x)Z^s+h\beta_n'(x))\right],\label{dersecyx} \\
\nonumber   \partial_{xx}&g(y,x)=\E\,\partial_{k+1k+1}f\big(y,x,\mathcal{E}^h_n(x,Z^s)\big)\\
  &+2\E\left[\partial_{k+1k+2}f(y,x,\mathcal{E}^h_n(x,Z^s))(1+\sqrt{h}\varsigma_n'(x)Z^s+h\beta_n'(x))\right]\notag\\&+\E\left[\partial_{k+2k+2}f(y,x,\mathcal{E}^h_n(x,Z^s))(1+\sqrt{h}\varsigma_n'(x)Z^s+h\beta_n'(x))^2+\partial_{k+2}f(y,x,\mathcal{E}^h_n(x,Z^s))\sqrt{h}\varsigma_n''(x)Z^s\right]\notag\\
  &+\E\left[\partial_{k+2}f(y,x,\mathcal{E}^h_n(x,Z^s))h\beta_n''(x)\right].\label{dersecx}\end{align}

 \smallskip
 \noindent {\sc Step 4}.
 For a Lipschitz continuous and directionally convex (non-decreasing in each of its variables) test function $f:\R^{k+2}\to\R$, the functions $f_m=\rho_m\star f$ are Lipschitz continuous and directionally convex (non-decreasing in each variable), infinitely differentiable with bounded derivatives and converge to $f$ uniformly as $m\to+\infty$, so that $\E \,f_m(y,x,\mathcal{E}^h_n(x,Z^s))\to \E \,f(y,x,\mathcal{E}^h_n(x,Z^s))$. Therefore it is enough to consider twice continuously differentiable test functions $f$ with bounded derivatives for which we may apply {\sc Step 3}. When $f$ is non decreasing in each of its variables, then clearly so is $g$ in each of its $k$ first variables. The function $g$ also is non-decreasing in its last variable by~\eqref{derxg} and~\eqref{derxcale}. When $f$ is directionally convex, so that its second order partial derivatives all are non-negative, then the first summands in the right-hand sides of~\eqref{dersecyy},~\eqref{dersecyx} and~\eqref{dersecx} are non-negative and the second summands in those of~\eqref{dersecyx} and~\eqref{dersecx} are non-negative by~\eqref{derxcale}. The fourth summand in the right-hand side of~\eqref{dersecx} is non-negative when $\beta$ is convex and $f$ is non-decreasing in its last variable. It is equal to $0$  disregarding the monotonicity of $f$ when $\beta$ is affine since $\beta_n$ is affine too and $\beta_n''$ vanishes. Therefore, to conclude that $g$ is directionally convex, it is enough to check that the third summand in the right-hand side of~\eqref{dersecx} is non-negative. To do so, we use the following equation based on an integration by parts for the second equality
\begin{align*}
   \E\left[\partial_{k+2}f(y,x,\mathcal{E}^h_n(x,Z^s))Z^s\right]&=\int_{-s}^{s}\partial_{k+2}f(y,x,\mathcal{E}^h_n(x,z))ze^{-\frac{z^2}{2}}\frac{dz}{\sqrt{2\pi}}\notag \\&=\int_{-s}^{s}\partial_{k+2k+2}f(y,x,\mathcal{E}^h_n(x,z))\sqrt{h}\varsigma_n(x)\left(e^{-\frac{z^2}{2}}-e^{-\frac{s^2}{2}}\right)\frac{dz}{\sqrt{2\pi}}
\notag\\&=\E\left[\partial_{k+2k+2}f(y,x,\mathcal{E}^h_n(x,Z^s))\sqrt{h}\,\varsigma_n(x){\mathbf 1}_{\{Z^s\ne 0\}}\Big(1-e^{\frac{(Z^s)^2-s^2}{2}}\Big)\right].\end{align*}
We deduce that the third summand is equal to 
\begin{align}
&\E\bigg[\partial_{k+2k+2}f(y,x,\mathcal{E}^h_n(x,Z^s))\!\!\left(\!\!(1+\sqrt{h}\,\varsigma_n'(x)Z^s+h\beta_n'(x))^2+h{\mathbf 1}_{\{Z^s\ne 0\}}\Big(1-e^{\frac{(Z^s)^2-s^2}{2}}\Big)\varsigma_n\varsigma_n''(x)\right)\!\bigg].\notag
\end{align}
When $c_\beta>0$, $\frac{\sqrt{c_\varsigma+2c_\beta}-\sqrt{c_\varsigma}}{2c_\beta}$ is the positive root of the quadratic function $\R\ni y\mapsto c_\beta y^2+\sqrt{c_\varsigma}y-\frac{1}{2}$. As a consequence, when $h\le\bar h_{\beta,\vartheta}$, $\frac{1}{2}-hc_\beta\ge \sqrt{hc_\varsigma }$ so that with $s\in\left[0,\frac{1}{2{\rm Lip}(\varsigma)\sqrt{h}}\right]$,~\eqref{minoderbetn},~\eqref{lipvarsig} and~\eqref{minosigsig''},
$$\forall x\in\R,\;(1+\sqrt{h}\,\varsigma_n'(x)Z^s+h\beta_n'(x))^2+h{\mathbf 1}_{\{Z^s\ne 0\}}\Big(1-e^{\frac{(Z^s)^2-s^2}{2}}\Big)\varsigma_n\varsigma_n''(x)\ge \left(\frac12-hc_\beta\right)^2-hc_\varsigma\ge 0.$$
This concludes the proof.\hfill $\Box$

\section{Proofs of the main results}\label{sec:proofs}
The next result is the key to get rid at the level of SDEs of the spatial semi-convexity of the square of one of the diffusion coefficients necessary for the comparison at the level of the truncated Euler schemes.

\begin{Proposition}\label{prop:approxeul}Assume that $b,\sigma,\beta,\vartheta$ satisfy~\eqref{coeflip} and $0\le\sigma\le \vartheta$. Then there is a sequence $(\sigma_n,\vartheta_n)_{n\ge 1}$ of approximations of $(\sigma,\vartheta)$ such that $\sigma_n,\vartheta_n:[0,T]\times\R\to\R$ are Lipschitz continuous in space uniformly in time and $n$,  $\mathbf{a}_{\sigma_n}\wedge\mathbf{a}_{\vartheta_n}>-\infty$ and $0\le \sigma_n\le\vartheta_n$ for each $n\ge 1$, and when $X_0,Y_0\in L^p(\P)$ for some $p\in[1,+\infty)$ the solutions $(X^n,Y^n)$ to
  $$X^n_t=X_0+\int_0^tb(s,X^n_s)ds+\int_0^t\sigma_n(s,X^n_s)dW_s\mbox{ and }Y^n_t=X_0+\int_0^t \beta(s,Y^n_s)ds+\int_0^t\vartheta_n(s,Y^n_s)dW_s,$$
  satisfy $\lim_{n\to+\infty}\Big\|\sup_{t\in[0,T]}|X_t-X^{n}_t|\Big\|_{p}=\lim_{n\to+\infty}\Big\|\sup_{t\in[0,T]}|Y_t-Y^{n}_t|\Big\|_{p}=0$.
\end{Proposition}

\noindent{\bf Proof of Proposition~\ref{prop:approxeul}.} Let $\rho$ be a $C^\infty$  probability density on the real line compactly supported on $[-1,1]$. We associate to $\rho$ its sequence of mollifiers  $\rho_n(x)=n\rho(nx)$, $n\in \N$. For $n\ge 1$ and $t\in[0,T]$, we set
\begin{equation*}
  \sigma_n(t,x)=\begin{cases}\rho_n\star\sigma(t,\cdot)(x)\mbox{ if } x\in[-n,n]\\
    \rho_n\star\sigma(t,\cdot)(n)+\left(0\vee\partial_x\rho_n\star\sigma(t,\cdot)(n)\right) (x-n)\mbox{ if }x\ge n\\
    \rho_n\star\sigma(t,\cdot)(-n)+\left(0\wedge \partial_x\rho_n\star\sigma(t,\cdot)(-n)\right)(x+n)\mbox{ if }x\le-n
  \end{cases},\end{equation*}
  and
  \begin{equation*}
  \vartheta_n(t,x)=\begin{cases}\rho_n\star\vartheta(t,\cdot)(x)\mbox{ if } x\in[-n,n]\\
    \rho_n\star\vartheta(t,\cdot)(n)+\left(0\vee\partial_x\rho_n\star\sigma(t,\cdot)(n)\vee\partial_x\rho_n\star\vartheta(t,\cdot)(n)\right) (x-n)\mbox{ if }x\ge n\\
    \rho_n\star\vartheta(t,\cdot)(-n)+\left(0\wedge \partial_x\rho_n\star\sigma(t,\cdot)(-n)\wedge \partial_x\rho_n\star\vartheta(t,\cdot)(-n)\right)(x+n)\mbox{ if }x\le-n
\end{cases}.
\end{equation*}
The inequality $0\le \sigma\le \vartheta$ clearly implies that $0\le \sigma_n\le \vartheta_n$ on $[0,T]\times[-n,n]$ and, with the above choice of slopes, $0\le \sigma_n\le \vartheta_n$ on $[0,T]\times\R$. For each $n\ge 1$ and $t\in[0,T]$, the function $x\mapsto \rho_n\star\sigma(t,\cdot)(x)$ (resp. $x\mapsto \rho_n\star\vartheta(t,\cdot)(x)$) is Lipschitz continuous with constant ${\rm Lip}(\sigma)$ (resp. ${\rm Lip}(\vartheta)$) so that its spatial derivative is bounded by ${\rm Lip}(\sigma)$ (resp. ${\rm Lip}(\vartheta)$) and $x\mapsto\sigma_n(t,x)$ (resp. $x\mapsto\vartheta_n(t,x)$) is Lipschitz continuous with constant ${\rm Lip}(\sigma)$ (resp. ${\rm Lip}(\sigma)\vee{\rm Lip}(\vartheta)$). Moreover since $\sigma$ satisfies~\eqref{coeflip} and $\rho_n$ and $\rho''_n$ are compactly supported and bounded,
$$\inf_{(t,x)\in[0,T]\times[-n,n]}\partial^2_{xx}(\rho_n\star\sigma(t,\cdot))^2(x)\ge 2\inf_{(t,x)\in[0,T]\times [-n,n]}(\rho_n''\star\sigma(t,\cdot))(\rho_n\star\sigma(t,\cdot))(x)>-\infty.$$
Since for each $t\in[0,T]$, the differences 
$$
2\rho_n\star\sigma(t,\cdot)(n)\left(\partial_x\rho_n\star\sigma(t,\cdot)(n)\right)^-\quad\mbox{ and }\quad 2\rho_n\star\sigma(t,\cdot)(-n)\left(\partial_x\rho_n\star\sigma(t,\cdot)(-n)\right)^+
$$ 
between  the right-hand slope and the left hand-slope of $x\mapsto \sigma_n^2(t,x)$ at point $x=n$  and $x=-n$ respectively,  is non-negative and $$\forall x\in (-\infty,-n)\cup(n,+\infty),\;\partial_{xx}\sigma^2_n(t,x)\in\left\{2\left(\partial_x\rho_n\star\sigma(t,\cdot)(n)\right)^2,0,2\left(\partial_x\rho_n\star\sigma(t,\cdot)(-n)\right)^2\right\},$$ we deduce that $\mathbf{a}_{\sigma_n}>-\infty$. We check in the same way that $\mathbf{a}_{\vartheta_n}>-\infty$.
Last, since $\rho$ is compactly supported on $[-1,1]$,
\begin{equation}\label{errsigman}
   \sup_{(t,x)\in[0,T]\times[-n,n]}|\sigma_n(t,x)-\sigma(t,x)|\le\frac{{\rm Lip}(\sigma)}{n}\mbox{ and }\sup_{(t,x)\in[0,T]\times[-n,n]}|\vartheta_n(t,x)-\vartheta(t,x)|\le\frac{{\rm Lip}(\vartheta)}{n}.
\end{equation}
Let us check that $\lim_{n\to+\infty}\Big\|\sup_{t\in[0,T]}|X_t-X^{n}_t|\Big\|_{p}=0$ when $X_0\in L^p(\P)$, the convergence of $Y^n$ to $Y$ following the same argument. We have 
$$X_t-X^n_t=\int_0^t(b(s,X_s)-b(s,X^n_s))ds+\int_0^t(\sigma_n(s,X_s)-\sigma_n(s,X^n_s))dW_s+\int_0^t(\sigma(s,X_s)-\sigma_n(s,X_s))dW_s.
$$
Thanks to the uniform in time and $n$ spatial Lipschitz continuity of the coefficients $(b,\sigma_n)$, standard estimations based on the H\"older and Burkholder-Davis-Gundy inequalities and Gronwall's lemma ensure for $p\ge 2$ the existence of a finite constant not depending on $m$ and $X_0\in L^p(\P)$ such that
\begin{align*}
  \Big \|\sup_{t\in [0,T]} |X_t-X^n_t|  \Big\|^p_p&\le C\int_0^T\|\sigma(s,X_s)-\sigma_n(s, X_s)\|^p_pds.\end{align*}
One has, using~\eqref{errsigman}
   and the affine growth of the coefficients $\sigma,\sigma_n$ in space uniform in time and $n\ge 1$,
\begin{align*}
   \sup_{t\in[0,T]}|\sigma(t,X_t)-\sigma_n(t,X_t)|\le \frac{{\rm Lip}(\sigma)}{n}+C\sup_{t\in[0,T]}|X_t|{\mathbf 1}_{\{\sup_{t\in[0,T]}|X_t|\ge n\}}.
\end{align*}
Reasoning like in the end of the proof of Proposition~\ref{propconvforteul2}, we deduce that for $p\ge 1$,
\begin{align*}
   \Big \|\sup_{t\in [0,T]} |X_t-X^n_t|  \Big\|_p\le C\left(\frac1n+\E^{1/p}\left[(1+|X_0|)^p\left(1\wedge \frac{1+|X_0|}{n}\right)^{\frac{1}{2}\wedge\frac{p}{4}}\right]\right)
\end{align*}
where the right-hand side goes to $0$ as $n\to+\infty$ since $X_0\in L^p(\P)$.
\hfill$\Box$

\bigskip
\noindent{\bf Proof of Theorem~\ref{thm:cvmarg}.} \noindent {\sc Step 1}. Let us temporarily assume either $(\mathcal X)$ or $(\mathcal Y)$ in place of $\forall t\in[0,T]$, $b(t,\cdot)$ is convex or $\forall t\in[0,T]$, $\beta(t,\cdot)$ is convex. We consider the truncated Euler schemes~\eqref{2eulertrunc} where for $k=0,\ldots,m-1$, $(\widetilde b^m_k(\cdot),\widetilde\sigma^m_k(\cdot),\widetilde\beta^m_k(\cdot),\widetilde\varsigma_k^m(\cdot))$ is equal to $$\left(\frac{m}{T}\int_{t_k}^{t_{k+1}}b(s,\cdot)ds,\sqrt{\frac{m}{T}\int_{t_k}^{t_{k+1}}\sigma^2(s,\cdot)ds},\frac{m}{T}\int_{t_k}^{t_{k+1}}\beta(s,\cdot)ds,\sqrt{\frac{m}{T}\int_{t_k}^{t_{k+1}}\vartheta^2(s,\cdot)ds}\right),$$
and with threshold $\mathbf{s}_m=\frac{\sqrt{m}}{2{\rm Lip}(\sigma)\sqrt{T}}$ when 
$(\mathcal X)$ holds or with threshold $\mathbf{s}_m=\frac{\sqrt{m}}{2{\rm Lip}(\vartheta)\sqrt{T}}$ when 
$(\mathcal Y)$ holds.
Since $b,\sigma,\beta,\vartheta$ satisfy~\eqref{coeflip} and $X_0$, $Y_0\!\in L^1(\P)$, it follows from~Proposition~\ref{prop:ConvfortEuler} that
\[
{\cal W}_1({\cal L}(X_{T}),{\cal L}(\widetilde X^m_{T}))+{\cal W}_1({\cal L}(Y_{T}),{\cal L}(\widetilde Y^m_{T}))\rightarrow 0\quad\mbox{as} \quad m\to+\infty.
\]
To establish the (increasing) convex ordering, we may restrict ourselves to Lipschitz (non-decreasing) convex functions $f:\R\to \R$ according to Proposition~\ref{prop:ConLip} 
. For such functions, by the dual formulation of the Wasserstein distance ${\cal W}_1$ with index $1$,     
\[
|\E\, f(\widetilde X^m_{_T})-\E \, f(X_{_T})|+ |\E\, f(\widetilde Y^m_{_T})-\E \, f(Y_{_T})|\le {\rm  Lip}(f) \left({\cal W}_1({\cal L}(X_{T}),{\cal L}(\widetilde X^m_{T}))+{\cal W}_1({\cal L}(Y_{T}),{\cal L}(\widetilde Y^m_{T}))\right)\to 0\]
as $m\to+\infty$.  So we can transfer  to the diffusion processes both convex orderings which hold for the truncated Euler schemes with $m$ large enough according to Proposition~\ref{prop:DirConvforEuler}. As for the  (non-decreasing) convexity of $x\mapsto \E\, f(X^x_{_T})$, it clearly holds for  any (non-decreasing) convex Lipschitz function $f$ owing to Remark~\ref{remdervar} and the above convergence. We can transfer the property to (non-decreasing) convex functions $f$ by following the reasoning in the last but one paragraph in the proof of Proposition~\ref{prop:DirConvforEuler}.

\medskip
\noindent {\sc Step 2}. According to {\sc Step 1}, the conclusion holds for $(X,Y)$ replaced by $(X^n,Y^n)$ given in Proposition~\ref{prop:approxeul} and it is enough to consider Lipschitz continuous test functions. We conclude since, by Proposition~\ref{prop:approxeul}, $\displaystyle \lim_{n\to+\infty}\Big\|\sup_{t\in[0,T]}|X_t-X^{n}_t|\Big\|_{1}=\lim_{n\to+\infty}\Big\|\sup_{t\in[0,T]}|Y_t-Y^{n}_t|\Big\|_{1}=0$. \hfill$\Box$

\medskip
\noindent{\bf Proof of Proposition~\ref{prop:CNconv2}.}
The second {and third statements are straightforward consequences} of the first one for the choice $X_0=\alpha x+(1-\alpha)y\le_{cvx}\varepsilon x+(1-\varepsilon) y=Y_0$ where $\varepsilon$ is distributed according to the Bernoulli law with parameter $\alpha$ and independent from $W$, for arbitrary $x,y\in\R$ and $\alpha\in[0,1]$.

Let us prove the first statement. For $s\in[0,T)$ and $t=s+h$ where $h=\frac Tm$ with $m\in{\mathbb N}$ such that $m\ge \frac{T}{T-s}$ so that $s+h\le T$, the convex ordering $(X_s,X_t)\preceq_{cvx} (Y_s,Y_t)$ implies that $\E[|X_{s+h}-X_s|]\le \E[|Y_{s+h}-Y_s|]$. It is enough to check that $\E[|X_{s+h}-X_s|]\sim\sqrt{\frac{2h}{\pi}}\E[|\sigma(s,X_s)|]$ and  $\E[|Y_{s+h}-Y_s|]\sim\sqrt{\frac{2h}{\pi}}\E[|\vartheta(s,Y_s)|]$ as $h\to 0$ to conclude, the case $s=T$ being recovered by continuity of $[0,T]\ni s\mapsto (\E[|\sigma(s,X_s)|],\E[|\vartheta(s,Y_s)|])$. Since $\sqrt{\frac{2h}{\pi}}\E[|\sigma(s,X_s)|]=\E\left[\left|\int_s^{s+h}\sigma(s,X_s)dW_u\right|\right]$, using the triangle inequality and Jensen's inequality, we get
\begin{align*}
  \left|\E\left[|X_{s+h}-X_s|-\sqrt{\frac{2h}{\pi}}|\sigma(s,X_s)|\right]\right|&\le \E\left[\left|\int_s^{s+h}b(u,X_u)du\right|\right]\\
  &\hskip 2.5cm  +\E^{1/2}\left[\left|\int_s^{s+h}(\sigma(u,X_u)-\sigma(s,X_s))dW_u\right|^2\right]\\
  &\le Ch(1+\|X_0\|_1)+\E^{1/2}\left[\int_s^{s+h}(\sigma(u,X_u)-\sigma(s,X_s))^2du\right].\end{align*}
Using standard estimations on the moments of the increments of $X$ under~\eqref{coeflip}, we obtain that
\begin{align*}\E\left[\int_s^{s+h}(\sigma(u,X_u)-\sigma(s,X_s))^2du\right]\le Ch^2(1+\|X_0\|_2)^2+2\E\left[\int_s^{s+h}(\sigma(u,X_s)-\sigma(s,X_s))^2du\right]
\end{align*}
and, when $X_0\in L^4(\P)$, reasoning like in the derivation of~\eqref{majodifftemps}, we check that
\begin{align*}
 \forall u\in[s,s+h],\;\|\sigma(u,X_s)-\sigma(s, X_{s})\|^2_2\le C\left(\varepsilon_{h,K}^2+  (1+\|X_0\|_4)^2\sqrt{1\wedge \frac{1+\|X_0\|_1}{K}}\right),
\end{align*}
where $K\ge 1$ and $\varepsilon_{h,K}:=\sup_{0\le s\le t\le (s+h)\wedge T,x\in[-K,K]}(|\sigma(s,x)-\sigma(t,x)|)$ goes to $0$ with $h$ by uniform continuity of $\sigma$ on the compact set $[0,T]\times[-K,K]$. Combining the three last inequalities with the subadditivity of the square root, we obtain that for deterministic initial conditions $x_0$,
\begin{align*}
   \left|\E\left[|X^{x_0}_{s+h}-X^{x_0}_s|-\sqrt{\frac{2h}{\pi}}|\sigma(s,X^{x_0}_s)|\right]\right|&\le C\sqrt{h}\left(\sqrt{h}(1+|x_0|)+\varepsilon_{h,K}+(1+|x_0|)\left(1\wedge \frac{1+|x_0|}{K}\right)^{\frac14}\right).
\end{align*}
Using Blagove$\check{\rm   s}\check{\rm  c}$enkii-Freidlin's theorem~\cite[Theorem 13.1, Section V.12-13, p.136]{RogersWilliamsII}, integrating with respect to the distribution of $X_0$ and using the independence between $X_0$ and $W$, we deduce that
\begin{align*}
   \left|\E\left[|X_{s+h}-X_s|-\sqrt{\frac{2h}{\pi}}|\sigma(s,X_s)|\right]\right|&\le C\sqrt{h}\left(\sqrt{h}(1+\|X_0\|_1)+\varepsilon_{h,K}+\E\left[(1+|X_0|)\left(1\wedge \frac{1+|X_0|}{K}\right)^{\frac14}\right]\right).
\end{align*}
For $X_0\in L^1(\P)$, by Lebesgue's theorem, the last expectation goes to $0$ as $K\to\infty$ so that we may choose $K$ such that it is arbitrarily small, while for fixed $K$ the sum of the first  two terms between brackets in the right-hand side goes to $0$ with $h$. Therefore the sum between these brackets goes to $0$ with $h$. 
The equivalent for $\E[|Y_{s+h}-Y_s|]$ is obtained in the same way.  \hfill$\Box$

\medskip
\noindent{\bf Proof of Theorem~\ref{thm:dirconv}.} The sequence $(X^n)_{n\ge 1}$ (resp. $(Y^n)_{n\ge 1}$) given by Proposition~\ref{prop:approxeul} converges in $L^p(\P)$ to $X$ (resp. $Y$). Since $F$ is continuous, we deduce that $(F(X^n))_{n\ge 1}$ (resp. $(F(Y^n))_{n\ge 1}$) converges in distribution to $F(X)$ (resp. $F(Y)$) and since $F$ has polynomial growth with order $p$, this sequence is uniformly integrable so that $\E \,F(X^n)\to\E \,F(X)$ (resp. $\E \,F(Y^n)\to\E \,F(Y)$) as $n\to+\infty$. Therefore we can suppose that either $(\mathcal X)$ or $(\mathcal Y)$ holds. We consider the truncated Euler schemes~\eqref{2eulertrunc} where for $k=0,\ldots,m-1$, $(\widetilde b^m_k(\cdot),\widetilde\sigma^m_k(\cdot),\widetilde\beta^m_k(\cdot),\widetilde\varsigma_k^m(\cdot))$ is equal to $$\left(\frac{m}{T}\int_{t_k}^{t_{k+1}}b(s,\cdot)ds,\sqrt{\frac{m}{T}\int_{t_k}^{t_{k+1}}\sigma^2(s,\cdot)ds},\frac{m}{T}\int_{t_k}^{t_{k+1}}\beta(s,\cdot)ds,\sqrt{\frac{m}{T}\int_{t_k}^{t_{k+1}}\vartheta^2(s,\cdot)ds}\right),$$
and with threshold $\mathbf{s}_m=\frac{\sqrt{m}}{2{\rm Lip}(\sigma)\sqrt{T}}$ when 
$(\mathcal X)$ holds or with threshold $\mathbf{s}_m=\frac{\sqrt{m}}{2{\rm Lip}(\vartheta)\sqrt{T}}$ when 
$(\mathcal Y)$ holds. Since $b,\sigma,\beta,\vartheta$ satisfy~\eqref{coeflip} and $X_0$, $Y_0\!\in L^p(\P)$, it follows from Corollary~\ref{coro:cvinterpol} that, as $m\to+\infty$,
  \[
 {\cal W}_p\left({\cal L}(X),{\cal L}\left(i_m\big((\widetilde X^m_{t_k})_{k=0,\ldots,m}\big)\right)\right)+{\cal W}_p\left({\cal L}(Y),{\cal L}\left(i_m\big((\widetilde Y^m_{t_k})_{k=0,\ldots,m}\big)\right)\right)\to 0,
\]
so that $\E \,F\big(i_m\big((\widetilde X^m_{t_k})_{k=0,\ldots,m}\big)\big)\to\E \,F(X)$ and $\E \,F\big(i_m\big((\widetilde Y^m_{t_k})_{k=0,\ldots,m}\big)\big)\to\E \,F(Y)$. Since the image of $\R_+^{m+1}$ by $i_m$ is included in $C([0,T],\R_+)$, $F\circ i_m$ is directionally convex on $\R^{m+1}$ and non-decreasing in each of its variables when $F$ is non-decreasing. With Proposition~\ref{prop:DirConvforEuler}, we deduce that for each $m\ge 1$, $\E \,F\big(i_m\big((\widetilde X^m_{t_k})_{k=0,\ldots,m}\big)\big)\le\E \,F\big(i_m\big((\widetilde Y^m_{t_k})_{k=0,\ldots,m}\big)\big)$. By taking the limit $m\to+\infty$ in this inequality, we conclude that $\E \,F(X)\le \E F(Y)$.

The (non-decreasing) convexity of $\R\ni x\mapsto\E \,F(X^x)$ when $F$ is continuous, (non-decreasing) directionally convex with polynomial growth is deduced by the argument in the first point of the remark which follows Theorem \ref{thmmarg1}.\hfill$\Box$

\small
\bibliographystyle{plain}\bibliography{propsdeconvex}

\appendix
\normalsize\section{Generalization of the stochastic calculus approach by El Karoui, Jeanblanc and Shreve in~\cite{EKJS}}\label{App:A} 
\begin{Proposition}\label{prop:NEKetal}
Let $b,\sigma: [0,T]\times \R\to \R$ satisfy \eqref{coeflip}. Assume furthermore that:

\medskip
\noindent $(i)$  the functions $b(t,\cdot)$ and $\sigma(t,\cdot)$ are continuously differentiable with derivatives $\partial_2b(t,\cdot)$ and $\partial_2\sigma(t,\cdot)$,  for every $t\!\in[0,T]$, 

\smallskip
\noindent $(ii)$ $\partial_2b(t,\cdot)$ and  $\partial_2\sigma(t,\cdot)$ are globally $\delta$-H\"older, $\delta\in (0,1]$, uniformly in $t\in[0,T]$,

\smallskip
\noindent $(iii)$ $\sigma\partial_2\sigma(t,\cdot)$ is Lipschitz uniformly in $t\in[0,T]$. 

\medskip
Let $X^x$ denote the solution of~\eqref{eq:SDE} starting from $x\!\in \R$. Let $f:\R^m\to\R$ be directionally convex (resp. directionally convex and non-decreasing in each of its variables) and let $0\le t_1\le t_2\le\cdots\le t_m\le T$,  $m\in\N$.

\medskip
\noindent $(a)$ Then,  $\R\ni x\mapsto \E f(X^x_{t_1},\cdots,X^x_{t_m})$ is convex (resp. convex and non-decreasing)
when $b(t,\cdot)$ is affine (resp. convex) for each $t\in[0,T]$.

\smallskip
\noindent $(b)$ In the particular case $m=1$, and $t_1=T$, this amounts to convexity (resp. non decreasing convexity) of $\R\ni x\mapsto \E f(X^x_{T})$ when $f:\R\to\R$ is convex (resp. non-decreasing and convex) and $b(t,\cdot)$ is affine (resp. convex) for each $t\in[0,T]$. 
\end{Proposition}

\noindent {\bf Remark.} Note that~\eqref{coeflip} implies that $\partial_2b$ and $\partial_2\sigma$ are bounded on $[0,T]\times\R$. 
{On the other hand, condition $(iii)$ implies that $\sigma^2(t,\cdot)$ has a Lipschitz derivative and hence is semi-convex.}

\bigskip
\noindent {\bf Proof.} {Since directional and regular convexity coincide for a function from $\R$ to $\R$, it is sufficient to prove $(a)$.} By Lemma \ref{lem:approxdirconv}, it is enough to deal with functions $f$ which are also Lipschitz continuous. By spatial convolution with $n^m\rho(nx)$ where $\rho$ is a compactly supported $C^\infty$ probability density on $\R^m$, such a function may be approximated uniformly by functions which are moreover $C^\infty$ and inherit its directional convexity and monotonicity properties. Therefore, we may even suppose that $f$ is $C^\infty$ with bounded first order derivatives.

By~\cite[Theorem 4.6.5]{Kunita} and its proof, there exists a modification of the solution to the stochastic differential equation~\eqref{eq:SDE} such that for $t\in[0,T]$, $x\mapsto X^x_t$ is a diffeomorphism of $\R$ with derivative $\partial_x X^x_t$ solving the equation obtained by spatial derivation of~\eqref{eq:SDE} :
$$d\partial_x X^x_t=\partial_2b(t,X^x_t)\partial_x X^x_t+\partial_2\sigma(t,X^x_t)\partial_x X^x_tdW_t,\;\partial_x X^x_0=1,$$
so that 
$$
\forall t\in[0,T],\;\partial_x X^x_t=e^{\int_0^{t}\partial_2 b(s,X^x_s)ds} e^{\int_0^{t}\partial_2 \sigma(s,X^x)dW_s -\frac12 \int_0^{t} (\partial_2 \sigma(s,X^x_s))^2ds}.
$$
 As a consequence $\partial_xf(X^x_{t_1},\cdots,X^x_{t_m})=\sum_{k=1}^m\partial_kf(X^x_{t_1},\cdots,X^x_{t_m})\partial_x X^x_{t_k}$. The boundedness of the first order derivatives of $f$ and of $\partial_2b$, $\partial_2\sigma$ implies that $\sup_{x\in\R}\E\left[\left(\sum_{k=1}^m\partial_kf(X^x_{t_1},\cdots,X^x_{t_m})\partial_x X^x_{t_k}\right)^2\right]<+\infty$ so that we may interchange the derivative and the expectation and obtain $$\partial_x\E\, f(X^x_{t_1},\cdots,X^x_{t_m})=\sum_{k=1}^m\E\Big(\partial_kf(X^x_{t_1},\cdots,X^x_{t_m})\partial_x X^x_{t_k}\Big).$$
By boundedness of $\partial_2\sigma$, $\E\,e^{\int_0^{t}\partial_2 \sigma(s,X^x)dW_s -\frac12 \int_0^{t} (\partial_2 \sigma(s,X^x_s))^2ds}=1$ for each $t\in[0,T]$. By Girsanov's theorem, under the probability measure  $\Q_{x,k}= e^{\int_0^{t_k}\partial_2 \sigma(s,X^x)dW_s -\frac12 \int_0^{t_k} (\partial_2 \sigma(s,X^x_s))^2ds}\cdot \P$, $\left(B^{x,k}_t = W_t -\int_0^t{\mathbf 1}_{\{s\le t_k\}} \partial_2 \sigma(s,X^x_s) ds\right)_{t\ge 0}$ is a Brownian motion and  $X^{x,k}$ is solution to the SDE
\[
X^{x,k}_t = x +\int_0^t (b+{\mathbf 1}_{\{s\le t_k\}}\sigma\partial_2\sigma)(s,X^{x,k}_s) ds +\int_0^t \sigma(s,X^{x,k}_s)dB^{x,k}_s.
\]Since $b,\sigma$ and $b+\sigma\partial_2\sigma$ are Lipschitz continuous in space uniformly in time,
${\cal L}_{\Q_{x,k}}(X^{x,k})= {\cal L}_{\P}(\tilde X^{x,k})$ where $\tilde X^{x,k}$ follows the same SDE as $X^{x,k}$ but with the Brownian motion $B^{x,k}$ replaced by $W$ :
\[
\tilde X^{x,k}_t = x +\int_0^t (b+{\mathbf 1}_{\{s\le t_k\}}\sigma\partial_2\sigma)(s,\tilde X^{x,k}_s) ds +\int_0^t \sigma(s,\tilde X^{x,k}_s)dW_s. 
\]
Then
\begin{align*}
   \partial_x\E f(X^x_{t_1},\cdots,X^x_{t_m})&=\sum_{k=1}^m\E_{\Q_{x,k}}\Big(\partial_kf(X^x_{t_1},\cdots,X^x_{t_m})e^{\int_0^{t_k}\partial_2 b(s,X^x_s)ds}\Big)\\&=\sum_{k=1}^m\E\Big(\partial_kf(\tilde X^{x,k}_{t_1},\cdots,\tilde X^{x,k}_{t_m})e^{\int_0^{t_k}\partial_2 b(s,\tilde X^{x,k}_s)ds}\Big).
\end{align*}
At this stage, since $b,\sigma$ and $b+\sigma\partial_2\sigma$ are Lipschitz continuous in space uniformly in time, by~\cite[Theorem 3.7 p.394]{RevuzYor}, $x\mapsto \tilde X^x$ is (a.s.) pathwise non-decreasing. For each $k\in\{1,\cdots,m\}$, since $\partial_{k\ell}f\ge 0$ for all $\ell\in\{1,\cdots,m\}$ by Property $\mathbf{P2}$ in Remark~\ref{remdircon}, $x\mapsto \partial_kf(\tilde X^{x,k}_{t_1},\cdots,\tilde X^{x,k}_{t_m})$ is non-decreasing as well. When $b(t,\cdot)$ is affine and therefore $\partial_2b(t,\cdot)$ constant for each $t\in[0,T]$, we deduce that $x\mapsto\partial_x\E f(X^x_{t_1},\cdots,X^x_{t_m})$ is non-decreasing and therefore $x\mapsto\E f(X^x_{t_1},\cdots,X^x_{t_m})$ is convex. The conclusion remains valid when $b(t,\cdot)$ is convex and therefore $\partial_2b(t,\cdot)$ non-decreasing for each $t\!\in [0,T]$ and $f$ is non-decreasing in each of its variables so that $\partial_k f\ge 0$ for each $k\in\{1,\cdots,m\}$. Finally, the pathwise non-decreasing property of $x\mapsto X^x$ deduced from~\cite[Theorem 3.7 p.394]{RevuzYor} implies that $x\mapsto\E f(X^x_{t_1},\cdots,X^x_{t_m})$ is non-decreasing when $f$ is non decreasing in each of its variables.\hfill$\Box$

\section{Proof of Example~\ref{exple:dirconv}$(d)$}\label{app:B}
The convexity of $\varphi$ implies that ${\cal C}([0,T], \R)\ni x\mapsto \int_0^T\varphi(x(s))ds$ is convex and the first statement follows from the convexity of the composition of a non-decreasing convex function with a convex function.
 To prove the second statement, let $x\in {\cal C}([0,T], \R)$ and $y, z\in {\cal C}([0,T], \R_+)$. For $s\in[0,T]$, we have $x(s)+y(s),x(s)+z(s)\in[x(s),x(s)+y(s)+z(s)]$ so that $\varphi(x(s)+y(s)),\varphi(x(s)+z(s))\in[\varphi(x(s)),\varphi(x(s)+y(s)+z(s))]$ by monotonicity of $\varphi$ and 
$\varphi(x(s)+y(s))+\varphi(x(s)+z(s))\le \varphi(x(s))+\varphi(x(s)+y(s)+z(s))$ by convexity of this function. As a consequence, $\int_0^T \varphi(x(s)+y(s))ds,\int_0^T \varphi(x(s)+z(s))ds\in\left[\int_0^T \varphi(x(s))ds,\int_0^T (\varphi(x(s)+y(s))+\varphi(x(s)+z(s))-\varphi(x(s)))\right]$ and 
$$
\int_0^T \varphi(x(s)+y(s))ds+\int_0^T \varphi(x(s)+z(s))ds-\int_0^T \varphi(x(s))ds\le \int_0^T \varphi(x(s)+y(s)+z(s))ds.
$$ 
With the convexity of $\Psi$ for the first inequality then its monotonicity for the second, we conclude that
\begin{align*}
   F(x+y)+F(x+z)&\le F(x)+\Psi\left(\int_0^T(\varphi(x(s)+y(s))+\varphi(x(s)+z(s))-\varphi(x(s)))ds\right)\\
   &\le F(x)+F(x+y+z).
\end{align*}

\end{document}